\pdfoutput=1

\documentclass[final]{siamltex}         % siam format but wider page
\usepackage[top=26mm, bottom=26mm, left=28mm, right=28mm]{geometry}
\usepackage{bm}        % bold math
\usepackage{color}    
\usepackage{soul}    

\usepackage{graphicx}				
\usepackage{amssymb}
\usepackage{amsmath}
\usepackage{multirow}
\usepackage{epstopdf}   % trying to auto generate PDF from EPS

\def\db{{\mathbf d}}%
\def\xb{{\mathbf x}}%
\def\kb{{\mathbf k}}%
\def\Fb{{\mathbf F}}%
\def\fb{{\mathbf f}}%
\def\sc{{\rm s_c} \, }%
\def\Cc{{\cal C}}%
\def\Mcb{\bm{\mathfrak M}}                        % all ref image
\def\Mcbh{{\widehat{\bm{\mathfrak M}}}}           % all ref image FT
\def\Mc{{\cal M}}                        % a ref image
\def\Ac{{\cal A}}                        % a ref image
\def\Mchat{{\widehat{{\cal M}}}}         % 2D FT of ref image
\newcommand{\kmax}{K}        % K max symbol
\newcommand{\bc}{\begin{center}}
\newcommand{\ec}{\end{center}}
\newtheorem{thm}{Theorem}
\newtheorem{rmk}{Remark}

\newcommand{\be}{\begin{equation}}
\newcommand{\ee}{\end{equation}}
\newcommand\bbR{\mathbb R}

\newcommand{\citep}[1]{\cite{#1}}      % since natbib incompat w/ siam bibsty

\begin{document}   %%%%%%%%%%%%%%%%%%%%%%%%%%%%%%%%%%%%%%%%%%%%%%%%%%%%%%

\title{Rapid solution of the cryo-EM reconstruction problem by 
frequency marching}

\author{Alex Barnett\thanks{Flatiron Institute, Simons Foundation, and 
Department of Mathematics, Dartmouth College}
\and
%{\em email:} {\tt ahb@math.dartmouth .edu}.}%
Leslie Greengard\thanks{Flatiron Institute, Simons Foundation, and
and Courant Institute, NYU}
\and
Andras Pataki\thanks{Flatiron Institute, Simons Foundation}
\and
Marina Spivak\thanks{Flatiron Institute, Simons Foundation}
}
\date{\today}	% Activate to display a given date or no date
\maketitle

\begin{abstract}
Determining the three-dimensional structure of proteins and protein complexes
at atomic resolution is a fundamental task in structural biology. Over the last
decade, remarkable progress has been made using ``single particle" 
cryo-electron microscopy (cryo-EM) for this purpose. 
In cryo-EM, hundreds of thousands of two-dimensional images 
are obtained of individual copies of the same particle, each held in a 
thin sheet of ice at some unknown orientation. Each image corresponds
to the noisy projection of the particle's electron-scattering density.
The reconstruction of a high-resolution image from this data is typically
formulated as a nonlinear, non-convex optimization problem for unknowns which 
encode the angular pose and lateral offset of each particle. 
Since there are hundreds of thousands of such parameters, this leads to a 
very CPU-intensive task---limiting both the number of particle images 
which can be processed and the number of independent reconstructions 
which can be carried out for the purpose of statistical validation.
Moreover, existing reconstruction methods typically require a good initial
guess to converge.
Here, we propose a deterministic method for high-resolution 
reconstruction that operates in an {\em ab initio} manner---that is, without
the need for an initial guess. It requires a predictable and relatively
modest amount of computational effort, 
by marching out radially in the Fourier domain from low to high frequency,
increasing the resolution by a fixed increment at each step.

\end{abstract}


\section{Introduction}

Cryo-electron microscopy (cryo-EM) is an extremely powerful technology for 
determining the three-dimensional structure of individual proteins and 
macro-molecular assemblies. Unlike X-ray diffraction-based methods, it does not
require the synthesis of high quality crystals, although there are still significant sample preparation issues involved. Instead, 
hundreds of thousands of copies of the same particle are held frozen in a thin sheet
of ice at unknown, random orientations. The sample is then placed 
in an electron microscope and, using a weak scattering approximation, the resulting
two-dimensional (2D) transmission images are interpreted as
noisy projections of the particle's electron-scattering density, with 
some optical aberrations that need to be corrected. 
Much of the recent advance in imaging quality has been as a result of hardware 
improvements (direct electron detectors), motion correction, and the
development of new image reconstruction algorithms.
Resolutions at the 2--4 \AA\ level are now routinely achieved.

There are a number of excellent overviews and texts in the literature 
to which we refer the reader for a more detailed introduction 
to the subject
\citep{Cheng2015,Frank2006,Milne2013,Nogales2016,VanHeel1997,Vinothkumar2016}.
One important issue to keep in mind, however, is that the field is currently lacking
a rigorous method for assessing the accuracy of the reconstructed density map,
although there are a variety of best practices in use intended to avoid over-fitting
\citep{Cheng2015,Rosenthal2015,Henderson2013,Heymann2015}.

In this paper, we concentrate on the reconstruction problem,
that is, converting a large set of noisy experimental cryo-EM 
images into a 3D
map of the electron-scattering density. 
As a rule, this is accomplished in two stages:
the generation of a low-resolution initial guess, 
either experimentally or computationally,
followed by an iterative
refinement step seeking to achieve the full available resolution permitted 
by the data. It is generally the refinement step which is CPU-intensive.
Although our {\em ab initio} method does not require an initial guess,
we briefly review some
of the existing approaches for the sake of completeness.

For methods which require a good, low-resolution initial guess
and the user has not supplied one experimentally,
many methods are based on the {\em common lines} principle.
This exploits the projection-slice theorem (see Theorem~\ref{fslicethm} and Fig.~\ref{figure:fourierslice}),
and was first introduced in \citep{vanHeel1987,Vainshtein1986}.
In the absence of noise, it is possible to show that any three randomly-oriented
experimental images can be used to define a coordinate system (up to mirror image
symmetry). Given that coordinate system, the angular assignment for all other images
can be easily computed. Unfortunately,
% because it uses only data along lines?
this scheme is 
very sensitive to noise. Significant improvement was achieved in the common lines
approach in \citep{Penczek1996}, where all
angles are assigned simultaneously by minimizing a global error. However, this
requires a very time-consuming calculation, as it involves searching in an
exponentially large parameter space of all possible orientations of
all projections. More recently, the SIMPLE algorithm \citep{Elmlund2012} 
used various optimization strategies, from simulated annealing to differential
evolution optimizers, to accelerate this search process.
Singer, Shkolnisky, and co-workers
\citep{Singer2009,Coifman2010,Singer2011,Shkolnisky2012}
proposed a variety
of rigorous methods based on the eigen-decomposition of certain sparse matrices 
derived from common line analysis and/or convex relaxation of the global 
least-squares error function.
Wang, Singer, and Wen
\citep{Wang2013} subsequently proposed a 
more robust self-consistency error, based on the sum of
absolute
residuals, which permits a convex relaxation and solution by semidefinite programming.

Another approach is the method of moments
\citep{Goncharov1987,Goncharov1988}, which is aimed at computing second order
moments of the density from second order moments of the 2D
experimental images. Around the same time,
Provencher and Vogel \citep{Provencher1988,Vogel1988}
proposed representing the electron-scattering density as a 
truncated expansion in orthonormal
basis functions in spherical coordinates and estimating the
maximum likelihood of this density. Recently, a machine-learning approach was
suggested, using a sum of
Gaussian ``pseudo-atoms" with unknown locations and radii
to represent the density, and Markov chain Monte Carlo sampling to 
estimate the model parameters \citep{Joubert2015}.

Once a low-resolution initial guess is established, 
the standard reconstruction packages
EMAN2 \citep{Bell2016,EMAN2}, SPIDER \citep{SPIDER}, SPARX \citep{SPARX},
RELION \citep{Scheres2012,Scheres2012b} and FREALIGN
\citep{Grigorieff2007,Grigorieff2016,Lyumkis2013} 
all use some variant of iterative projection matching in either physical
or Fourier space.
Some, like RELION, SPARX and EMAN2 make use of soft matching or a 
regularized version of the maximum likelihood 
estimation (MLE) 
framework
introduced by Sigworth \citep{Sigworth1998,Sigworth2010}.
For this, a probability distribution for angular assignments is determined
for each experimental image.
Others, like FREALIGN, assign unique angular assignments to each image.
Two drawbacks of these methods are that they can be time-consuming, 
especially those based on MLE, and it can be difficult to determine if
and when a global optimum has been reached.

In order to accelerate the MLE-based methods,
it was suggested in \citep{Dvornek2015} that particle images and 
structure projections be represented in low-dimensional
subspaces that permit rotation, translation and comparison by
defining suitable operations on the subspace
bases themselves. They demonstrated 300-fold speedups in reconstruction.
Recently,
Brubaker, Punjani, and Fleet
\citep{Brubaker2015a} introduced a new scheme
based on a probabilistic generative model, marginalization over angular
assignments, and optimization using stochastic gradient descent and importance sampling.
They demonstrated that their method was both efficient and insensitive to the 
initial guess.
Finally, there has been significant effort aimed at harnessing high-performance
computing hardware including GPUs for the most compute-intensive tasks in the 
cryo-EM reconstruction pipeline \citep{Kimanius2016}.

In this paper, we will focus on a new method for refinement, assuming
that all particle images are drawn from a homogeneous population. 
One of our goals is the creation of a method for refinement that is sufficiently
fast that it can be run multiple times on the same data,
opening up classical jackknife, bootstrap, or cross-validation
statistics to be used for validation and resolution assessment.
At present, the gold standard is based on methods such as
``Fourier shell correlation''
\citep{Henderson2012,Penczek2010rev,Rosenthal2015},
which can be viewed as a jackknife with two samples.

The basic intuition underlying our approach is already shared with many of the
standard software packages, such as RELION and FREALIGN, as well as the stochastic
optimization method of \citep{Brubaker2015a}: namely, that 
the best path to a refined structure is achieved by
{\em gradually increasing resolution}.
The main purpose of the present paper is to propose a deterministic and mathematically
precise version of this idea, carried out in the Fourier domain.
Unlike existing schemes, it involves no global optimization.
Instead, we use resolution (defined by the maximum frequency content of the 
current reconstruction) as a homotopy path. For each small
step along that path, we solve only {\em uncoupled} projection matching problems
for each experimental image. 
More precisely, let $k \in [0,K]$ denote the band-limit of the
model at the current step of the reconstruction algorithm, where $K$ denotes the 
maximum resolution we seek to achieve.
For the model, we generate a large number of ``templates''---
simulated projections of the model---at a large number of orientations.
For each experimental image, we then find the template that is a best match
and assign the orientation of that template to the image. Given the current
angular assignments of all experimental images, we solve a linear least squares
problem to build a {\em new} model at resolution $k+ \Delta k$ and 
repeat the process until the maximum resolution $K$ is reached
(see section \ref{RLsec}).

With $M$ images, each at a resolution of $K \times K$ pixels, our scheme
requires $O(M K^4)$ or
$O(M K^5)$ work, depending on the cost of template-matching. If this is 
done by brute force, the second estimate applies. If a hierarchical 
but local search strategy is employed for template matching, then the
cost of the least squares procedures dominates and the $O(M K^4)$ complexity
is achieved. 
(The memory requirements are approximately $8M K^2 + 8K^3$ bytes.)

\begin{rmk}
Some structures of interest have non-trivial point group
symmetries. As a result, there may not be a unique angular assignment
for each experimental image. In the most extreme case, one could 
imagine imaging perfect spheres, for which angular assignment makes no sense
at all. Preliminary experiments with noisy projection data have been successful
in this case, using the randomized assignment scheme discussed below 
in section \ref{s:match}. We believe that our procedure is easily modified 
to handle point group symmetries as well but have not explored this class
of problems in detail (see \cite{Belnap2010} for further discussion).
\end{rmk} 

Our frequency marching scheme
was inspired by the method of recursive
linearization for acoustic inverse scattering, originally introduced by Y. Chen 
\citep{BaoLi2015,Chenrep1088,Chen}.
We show here that high resolution and low errors can be achieved systematically, with a well-defined
estimate of the total work required.
As noted above, unlike packages such as RELION and EMAN2, we do not address 
the issues that arise when there is heterogeneity in the data sets due to the 
presence of multiple quasi-stable conformations of the particles being imaged.
We also assume that the experimental images have known in-plane translations.
Although fitting for such translations would be desired in a production code,
we believe that this can be included in our algorithm
with only a small constant factor increase in computation time
{(see the concluding section for further discussion).}
We do, however, include in our algorithm 
a known {\em contrast transfer function} (CTF) which models
realistic aberrations of each experimental image.

In sections \ref{prelimsec}--\ref{lsqsec}, we 
introduce the notation necessary to describe the algorithm in detail, as well
as the various computational kernels that will be needed. The
frequency marching
procedure is described in section \ref{RLsec}.
Section \ref{resultsec1}
presents our numerical results,
which use simulated data derived from known atomic positions for
three relevant protein geometries of interest.
The use of simulated data allows us to investigate the effect of
signal-to-noise ratio (SNR) on reconstruction quality.
We draw conclusions in section~\ref{conc}.

\section{Mathematical preliminaries} \label{prelimsec}

We begin by establishing some notation. Throughout this paper, the 
unknown electron-scattering density will be denoted by $f(\xb)$ where
$\xb = (x,y,z)$ in Cartesian coordinates.
We assume, without loss of generality, that the unit of length is chosen
so that the particle (support of $f$) fits in the unit ball at the origin. 
The Fourier transform of $f(\xb)$ will be 
denoted by $F(\kb)$, where
$\kb = (k_1,k_2,k_3)$ in Cartesian coordinates and 
$\kb = (k,\theta,\phi)$ in spherical coordinates, with 
$k_1 = k \sin \theta \cos \phi$, $k_2 = k \sin \theta \sin \phi$, 
and $k_3 = k \cos \theta$.
Since we will use a spherical discretization of Fourier space, we write 
the standard Fourier relations in the form:
\begin{equation}
F(k,\theta,\phi) = 
\int_{-1}^1
\int_{-1}^1
\int_{-1}^1
f(x,y,z) e^{i k ( x \sin \theta \cos \phi + \,
y \sin \theta \sin \phi  + \,
z \cos \theta )} \, dx \, dy \, dz
\label{FTdef}
\end{equation}
and
\begin{equation}
f(x,y,z) = 
\frac{1}{(2\pi)^3}
\int_{0}^{\infty}
\int_{0}^{2 \pi}
\int_{0}^{\pi}
F(k,\theta, \phi)  e^{- i k (x \sin \theta \cos \phi  + \,
y \sin \theta \sin \phi  + \,
z \cos \theta )} \, k^2 \sin \theta \, d\theta \, d\phi \, dk.
\label{FTinvdef}
\end{equation}
We will make use of Clenshaw-Curtis (Chebyshev) quadrature
in $\cos \theta$ for the inner integral
(see section \ref{sec:disc}),
since this is spectrally accurate for smooth integrands and
results in
nodes that are equispaced in the parameter $\theta$, which simplifies
the task of interpolation.

\subsection{Rotation and projection operators}

Let $(1,\alpha,\beta)$ denote the spherical coordinates of a point on the unit 
sphere, which we will also view as an {\em orientation vector}. 
With a slight abuse of notation, we 
will often identify $(\alpha,\beta)$ with the vector 
$(1,\alpha,\beta)$. 
Rather than assuming the electron beam orientation is along the $z$-axis and that the
particle orientations are unknown,
it is convenient to imagine that the particle of interest is fixed in the laboratory 
frame and that each projection obtained from electron microscopy corresponds
to an electron beam in some direction $(\alpha,\beta)$.

\begin{definition}  
Let $\db = (1,\alpha,\beta)$ be an orientation vector, and let
$(r,\psi,s)$ denote a 
cylindrical coordinate system in $\bbR^3$, where $(r,\psi)$ are 
polar coordinates in the projection plane orthogonal to $\db$ and
$s$ is the component along $\db$.
The projection of the function $f(\xb)$ in the direction 
$\db$ is denoted by $P_{\alpha,\beta}[f]$ with
\[ P_{\alpha,\beta}[f](r,\psi) = \int_{-\infty}^\infty 
f(r,\psi,s) \, ds \, . \]
\end{definition}

\begin{figure}
\begin{center}
\includegraphics[width=3.5in]{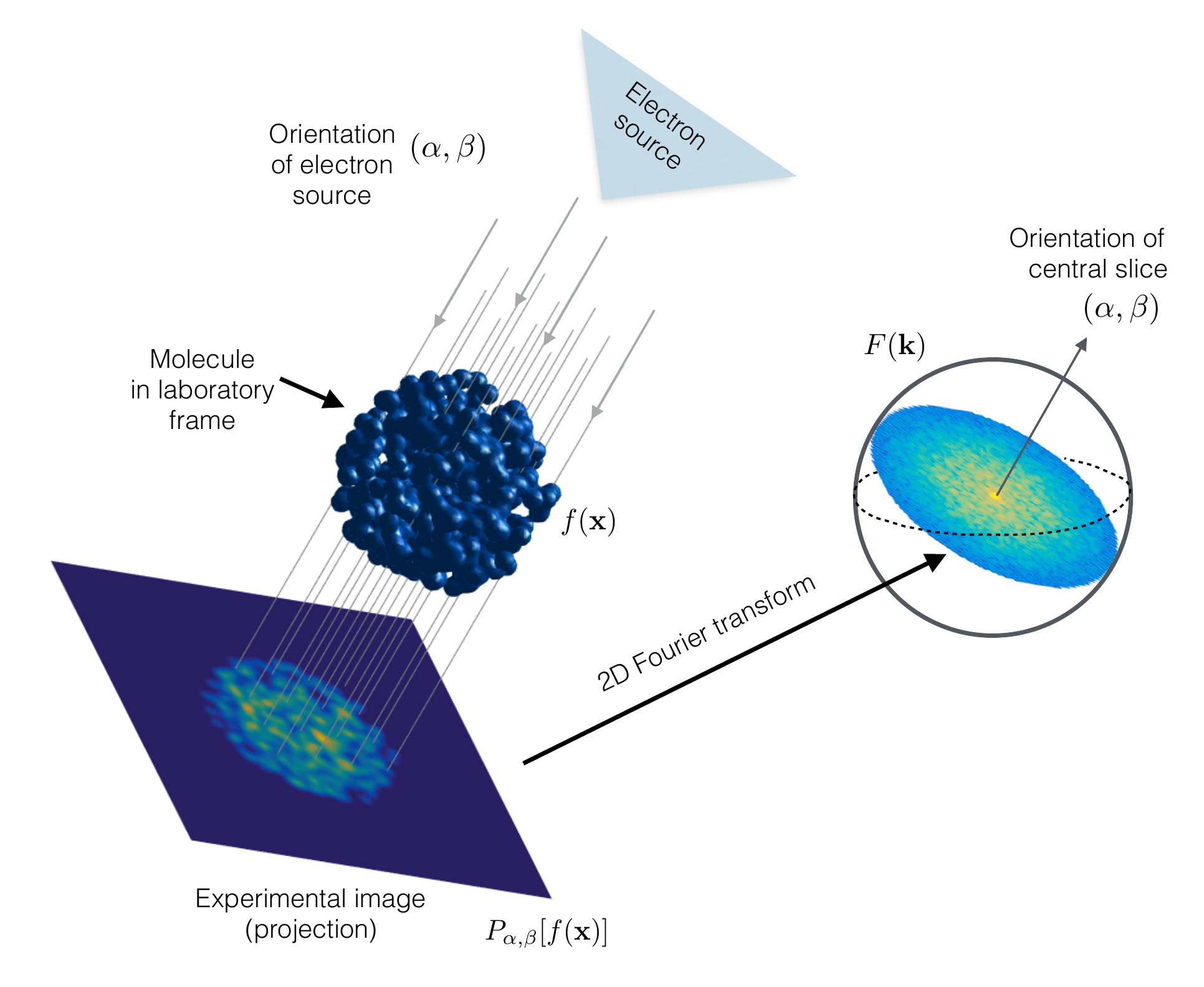}
\end{center}
\caption{Illustration of the projection-slice theorem:
the 2D Fourier transform of the projection of a compactly 
supported function $f(\xb)$ equals a central slice of the Fourier transform
$F(\kb)$ through the origin, with normal vector $(1,\alpha,\beta)$.}
\label{figure:fourierslice}
\end{figure}

There is a simple connection between the projection of a function $f$,
namely $P_{\alpha,\beta}[f]$, and its Fourier transform $F$, given by 
Theorem \ref{fslicethm} and illustrated in Fig.~\ref{figure:fourierslice}. 
For this, we will need to define a central slice of $F$.

\begin{definition}  
Let $\db = (1,\alpha,\beta)$ be an orientation vector and
let $F(\kb)$ denote a function in the three-dimensional Fourier
transform domain. Then the {\em central slice}
$S_{\alpha,\beta}[F]$ is defined to be the restriction of $F$ to the plane through 
the origin with normal vector $\db$.
\end{definition}

\begin{thm} \label{fslicethm} {\em (The projection-slice theorem)}.
Let $f(\xb)$ denote a compactly supported function in $\bbR^3$, let 
$F(\kb)$ denote its Fourier transform and let 
$(\alpha,\beta)$ denote an orientation vector.
Let $P_{\alpha,\beta}[f]$ denote the corresponding projection of $f$ 
and let $\widehat{P_{\alpha,\beta}[f]}$ denote its 2D
Fourier transform (with normalization analogous to \eqref{FTdef}--\eqref{FTinvdef}.
Then $\widehat{P_{\alpha,\beta}[f]}$ corresponds to an equatorial
slice through $F(\kb)$, with orientation vector $(\alpha,\beta)$. That is, 
\[ \widehat{P_{\alpha,\beta}[f]} = S_{\alpha,\beta}[F]. \]
\end{thm}

The proof of this theorem is straightforward (see, for example, \citep{Natterer}).
There is a third angular degree of freedom which must be taken into 
account in cryo-EM, easily understood 
by inspection of Fig.~\ref{figure:fourierslice}.
In particular, any rotation of the image plane about the orientation
vector $(\alpha,\beta)$ could be captured during the experiment. 

To take this rotation into account, we need to be more precise
about the polar coordinate system $(r,\psi)$ in the image plane. We will
denote by $\psi$ the angle subtended in the image plane
relative to the projection of the negative $z$-axis
in the fixed laboratory frame.
(In the case of projections with directions passing through the poles
$\alpha = 0$ and $\pi$, the above
definition becomes ambiguous, and one may set $\psi$ to be
the standard polar angle in the $xy$ plane.)
Then $\gamma$ specifies in-plane rotation of a projection.
Thus, the full specification of an arbitrary projection takes the form
\be
P_{\alpha,\beta,\gamma}[f](r,\psi) \; = \; P_{\alpha,\beta}[f](r,\psi-\gamma)
~.
\label{Pabg}
\ee

\begin{rmk}
The angles $(\alpha,\beta,\gamma)$ correspond to a particular 
choice of Euler angles that define an arbitrary rotation of a rigid body 
in three dimensions. The notation introduced here is most convenient for the 
purposes of our method {(corresponding to an extrinsic rotation
of $\alpha$ about the $y$ axis, an extrinsic rotation of $\beta$ about the 
$z$ axis, and an intrinsic rotation of $\gamma$ about the 
new $z$ axis).}
\end{rmk}

We will denote the set of experimental images obtained from electron microscopy
(the input data) by $\Mcb$, with 
the total number of images given by $M = |\Mcb|$.
Each image $\Mc^{(m)} \in \Mcb$ has support in the unit disc.
We will often omit the superscript $(j)$ when the context is clear.
Representing the image in Cartesian coordinates,
we let $\Mchat^{(m)}$ denote its 2D Fourier 
transform in polar coordinates:
\begin{equation}
\Mchat^{(m)}(k,\psi) =
\int_{-1}^1
\int_{-1}^1
\Mc^{(m)}(x,y) e^{i k (\cos \psi \, x + \, \sin \psi \, y) } \, dx \, dy.
\label{FT2def}
\end{equation}

% Dddddddddddddddddddddddddddddddddddddddddddddddddddddddddddddddd
\subsection{Discretization} \label{sec:disc}

We discretize the full three-dimensional
Fourier transform domain ($\kb$-space) in spherical 
coordinates as follows.
We choose $N_r$ equispaced quadrature nodes
in the radial direction, between zero
and a maximum frequency $\kmax$. The latter sets the achievable resolution.
On each sphere of radius $k$, we form a 2D product grid from
$N_\phi$ equispaced points in the $\phi$ direction, and
the following $N_\theta$ equispaced points in the $\theta$ direction:
\[  \theta_j = \frac{(2j-1) \pi}{2N_\theta} ~, \qquad 1 \leq j \leq N_\theta ~.
\]
Writing the volume element $k^2 \sin \theta \,d\theta\, d\phi\, dk$ as
$k^2 d\mu \,d\phi\, dk$ where $\mu = \cos \theta$ is the scaled $z$-coordinate,
we note that the scaled $z$-coordinates of the nodes $\mu_j = \cos \theta_j$
are thus located at the classical (first kind) Chebyshev nodes
on $[-1,1]$.
The advantage of this particular choice for $\theta$ (or $\mu$)
nodes is that they will be convenient later for local
interpolation in $\theta$ near the poles.

The complete spherical product grid of $N_r N_\theta N_\phi$ nodes
is used for $\kb$-space quadrature,
with the weight associated with each node being $(2\pi/N_\phi) k^2 w_j$, where
$w_j$ is the weight corresponding to the node $\mu_j$ on $[-1,1]$.
(Here the $\mu$ quadrature is sometimes known as Fej\'er's first rule
\cite{kink}.)
In particular, this quadrature will be used in computing the
final density in physical space by means of the inverse Fourier transform
\eqref{FTinvdef}.

It is important to note that, because the support of $f$ is in the unit ball,
we have very precise bounds on the {\em smoothness} of $F(\kb)$.
Namely, $F$, as a function of $\kb$, is bandlimited to unit ``frequency''
(note that here, and here only, we use ``frequency'' in the reverse sense
to indicate the rate of
oscillation of $F$ with respect to the Fourier variable $\kb$).
Furthermore, for targets $\xb$ in the unit ball, the same is true of
the exponential function in \eqref{FTinvdef}.
Thus, the integrand is bandlimited to a ``frequency'' of 2.
This means that the above quadrature scheme is
superalgebraically convergent with respect to $N_\phi$ and $N_\theta$,
due to well-known
results on the periodic trapezoid and Chebyshev-type quadratures \cite{kink}.
Furthermore, for oscillatory periodic band-limited functions, the
periodic trapezoid rule reaches full accuracy once
``one point per wavelength'' is exceeded for the integrand \cite{PTRtref}
(note that this is half the usual Nyquist criterion).
These results are expected to
carry over to the spherical sections of 3D bandlimited functions that we
use, allowing for superalgebraic error terms.
Since the most oscillation with respect to $\theta$ and $\phi$
occurs on the largest $k$ sphere,
the above sampling considerations imply
$N_\phi \ge 2 \kmax$ and $N_\theta \ge \kmax$.
In practice, to ensure sufficient accuracy,
we choose values slightly (i.e., up to a factor 1.5) larger than these bounds.

\begin{rmk}
In practice, the sphere can be sampled more uniformly and more efficiently
by choosing $N_\theta$ to vary with $k$ and by reducing the number of 
azimuthal points $N_\phi$ near the poles. For simplicity of presentation,
we keep them fixed here (see section \ref{conc} for further discussion).
\end{rmk}

Quadrature in the $k$ direction is only second order accurate,
as we are relying on the trapezoidal rule with an equispaced grid.
Spectral convergence could be achieved using Gauss-Legendre quadrature,
but we find that accuracy with a regular grid is sufficient
using a node spacing $\delta k := \kmax/N_r \approx 2$,
so that $N_r=O(\kmax)$.

In our model, rather than using the spherical grid points themselves to sample
$F(\kb)$, we will, for the most part, make use of a spherical harmonic
representation.
That is, for each fixed radial value $k$, we will represent 
$F(k,\theta,\phi)$ on the corresponding spherical shell in the form
\begin{equation}
\label{shellrep}
F(k,\theta,\phi)  = 
 \sum_{n=0}^{p(k)} \sum_{m=-n}^n f_{nm}(k) Y_n^m(\theta,\phi) \, .
\end{equation}
Here,
\begin{equation}
Y_n^m(\theta,\phi) =
S_n^m P_n^{|m|}(\cos \theta)
		      e^{i m \phi} \, ,
\label{ynmpnm}
\end{equation}
where $P_n(x)$ denote the standard Legendre polynomial of degree $n$,
the associated Legendre functions $P_n^m$ are defined by the Rodrigues' formula
\[ P_n^m(x) = (-1)^m (1-x^2)^{m/2} \frac{d^m}{dx^m} P_n(x), \] 
and
\begin{equation}
S_n^m =
\sqrt{\frac{2n+1}{4 \pi}}
\sqrt{\frac{(n-|m|)!}{(n+|m|)!}} \, .
\label{snm}
\end{equation}

The functions $Y_n^m$ are orthonormal with respect to the $L_2$ inner
product on the unit sphere
\begin{equation}\label{bil_s}
\langle f,g \rangle =\int_0^{2\pi} \int_0^{\pi} f(\theta,\phi) g(\theta,\phi) \sin \theta \, d\theta\, d\phi \, .
\end{equation}

\begin{rmk}
Note that the degree of the expansion in \eqref{shellrep} is a function of 
$k$, denoted by $p(k)$. 
It is straightforward to show that spherical 
harmonic modes of degree greater than $k$ are exponentially decaying on a
sphere of radius $k$ given that the original function $f(x,y,z)$ is supported
in the unit ball in physical space. 
This is the same argument used in discussing {\em grid} resolution above
for $N_r$, $N_\phi$ and $N_\theta$.
In short, an expansion of degree $O(k)$ is sufficient and, throughout this
paper, we simply fix the degree $p(k) = k+2$.
\end{rmk}

\begin{rmk}
The orthogonality of the $Y_n^m$ allows us to obtain the coefficients of the 
spherical harmonic expansion $f_{nm}$ above via projection. Moreover, 
separation of variables permits all the necessary integrals to be computed
in $O(p(k)^3)$ operations.
\end{rmk}

\section{Template generation}

\begin{figure} 
\begin{center}
\includegraphics[width=3in]{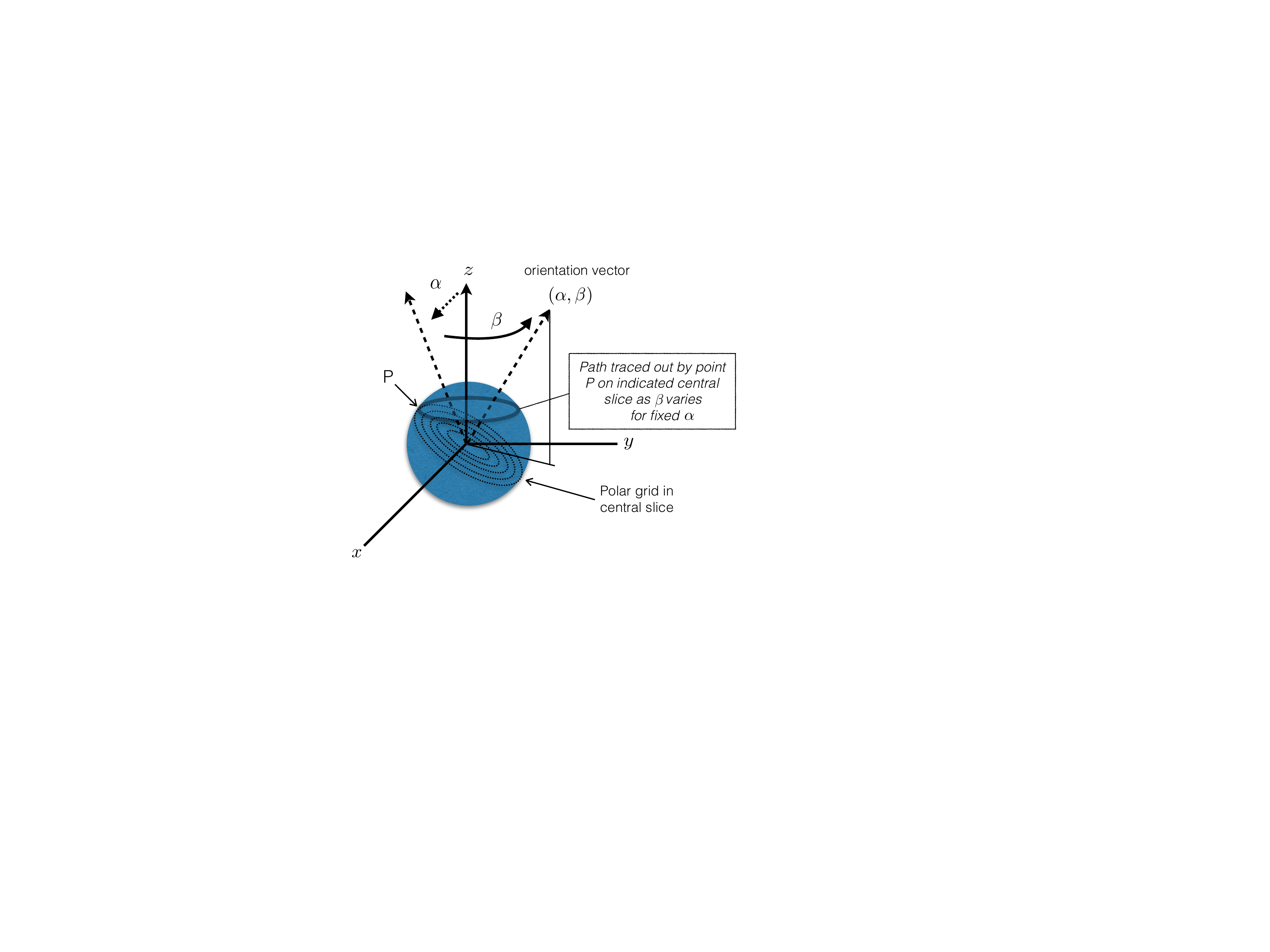}
\end{center}
\caption{To generate templates for projection matching,
we need to sample the central slice $S_{\alpha,\beta}[F]$ on a polar grid
for every orientation vector $(\alpha,\beta)$.
The process is accelerated by observing that the $z$-coordinate
for any such target point depends only on $\alpha$ and not on $\beta$.
This permits the effective use of separation of variables
(eq. \eqref{compall}).}
\label{figure:templategen}
\end{figure}

In any projection matching procedure, a recurring task is that of 
{\em template generation}. That is, given a current model 
defined by $f(\xb)$ or $F(\kb)$, we must generate a collection of projection images 
of the model for a variety of orientation vectors. We will denote those
orientations by $(\alpha_i,\beta_j)$, leaving $(\theta,\phi)$ to refer to the 
angular coordinates in Fourier space for a fixed frame of reference.
From the projection-slice theorem (Theorem \ref{fslicethm}), we have
\[
\widehat{P_{\alpha_i,\beta_j}[f]} = S_{\alpha_i,\beta_j}[F]~.
\]
In fact, our projection matching procedure will work entirely
in 2D Fourier space, ie
using $\widehat{P_{\alpha_i,\beta_j}[f]}(k,\psi)$,
where $(k,\psi)$ are polar coordinates
in the projection plane,
so our task is
simply to compute the central slice $S_{\alpha_i,\beta_j}[F]$ on a polar grid
for each $i$ and $j$ (see Fig. \ref{figure:templategen}).
The real-space projections $P_{\alpha_i,\beta_j}[f]$ will not be needed.

In order to generate all templates efficiently,
let a discrete set of desired orientation vectors be denoted by 
$\{ (\alpha_i,\beta_j)|\, i = 1,\dots,N_\theta, j = 1,\dots,N_\phi \}$.
We have found that, to achieve acceptable errors at
realistic noise levels, these grid sizes $N_\theta$ and $N_\phi$
can be chosen to be the same as for the quadrature points via
the Nyquist criterion in the previous section.
Each $\alpha$ corresponds to an initial rotation about the $y$-axis in the 
laboratory frame, while each $\beta$ corresponds to a subsequent rotation
about the $z$-axis in the laboratory frame.
(The fact that projection matching can be simplified by creating a sufficiently
fine grid on the sphere was already discussed and used in 
\citep{Penczek1992,Penczek1994}.)

To generate the points at a given radius $k$ within each slice,
we must parametrize points in a polar coordinate system
on the equatorial plane in the laboratory frame that have been rotated
through the preceding actions. 
Using Cartesian coordinates for the moment, it is easy to see that
{\small 
\[ \left( \begin{array}{c} 
k \cos \psi \\
k \sin \psi \\
0 
\end{array}
\right)  \ \xrightarrow{\alpha\ {\rm rotation}} \ 
\left( \begin{array}{c} 
k \cos \alpha \cos \psi \\
k \sin \psi \\
-k \sin \alpha \cos \psi 
\end{array}
\right)  \ \xrightarrow{\beta\ {\rm rotation}} \ 
\left( \begin{array}{c} 
k[ \cos \beta \cos \alpha \cos \psi - \sin \beta \sin \psi] \\ 
k [\sin \beta \cos \alpha \cos \psi + \cos \beta \sin \psi] \\ 
k [- \sin \alpha \cos \psi] 
\end{array} \right).
\]
}
Note that the $z$ coordinate (and therefore the 
spherical coordinate $\theta$) in the laboratory frame is 
independent of $\beta$, while
the $x,y$ coordinates in the laboratory frame depend on $\alpha$, $\beta$ and 
$\psi$.

For each sampled normal orientation vector $(\alpha_i,\beta_j)$ for the slice,
let 
$\{ \psi_l | l = 1,\ldots,N_\phi\}$
denote equispaced points on $[0,2\pi]$.
From above, the sample points with radius $k$,
which are located at $(k \cos \psi_l,k \sin \psi_l,0)$
before rotation, move to
\[ 
\left( \begin{array}{c} 
x_{ijl} \\ y_{ijl} \\ z_{il}
\end{array} \right)
= 
\left( \begin{array}{c} 
k [\cos \beta_j \cos \alpha_i \cos \psi_l - \sin \beta_j \sin \psi_l] \\ 
k [\sin \beta_j \cos \alpha_i \cos \psi_l + \cos \beta_j \sin \psi_l] \\ 
k [- \sin \alpha_i \cos \psi_l] 
\end{array} \right).
\]

Fixing the radius $k$ for now,
we denote the spherical coordinates of these points 
$(x_{ijl}, y_{ijl},z_{il})$ by 
$(k, \theta_{il}, \phi_{ijl})$.
To generate the template data, we now need to evaluate
the spherical harmonic expansion at these points, namely
\begin{align}\label{compall}
F(k,\theta_{il},\phi_{ijl}) 
&= 
 \sum_{n=0}^{p(k)} \sum_{m=-n}^n f_{nm}(k) Y_n^m(\theta_{il},\phi_{ijl}) 
\nonumber \\
&= 
 \sum_{m=-n}^n A_{iml} e^{i m \phi_{ijl}} 
\end{align}
for $i = 1,\dots,N_\theta$ and $j = 1,\dots,N_\phi$,
where
\[ A_{iml} = 
 \sum_{n=|m|}^{p(k)}f_{nm}(k) S_n^m P_n^m(\cos \theta_{il}). 
\]
Note that only the degrees at least as large as the magnitude of the
order $|m|$ contribute.
For each radius $k$, the cost of computing the set $\{ A_{iml} \}$ is 
$O(N_\theta N_\phi p(k)^2)$ and the cost of 
the subsequent evaluation of all values 
$F(k,\theta_{il},\phi_{ijl})$ is  
$O( N_\theta N_\phi^2 p(k))$.
Since $p(k)= O(\kmax)$,
both terms have the 
asymptotic complexity $O(\kmax^4)$. This generates $O(\kmax^3)$
template data points on a single spherical shell.
Then, summing over all $k$ shells,
the cost to generate the full spherical grid templates up to 
a resolution of $\kmax$ requires $O(\kmax^5)$ work, and generates
$O(\kmax^4)$ data.

Letting $k_q \in \{ k_1,\dots,k_{N_r} \}$
denote the discretization in $k$, the above procedure
evaluates the $O(\kmax^2)$ samples comprising the polar grid $(k_q, \psi_l)$,
for each of the $O(\kmax^2)$ central slices.
However, it is more convenient to store these templates in
terms of their angular Fourier series for each slice,
ie,
\be
S_{\alpha_i,\beta_j}[F](k_q,\psi) = 
\sum_{n=-N_q}^{N_q} S^{ij}_n(k_q) e^{i n \psi},
\label{Sijnk}
\ee
with $2N_q +1$ Fourier modes used on the $q$th ring.
By band-limit considerations of the spherical harmonics,
one need only choose $N_q = p(k_q)$, the maximum degree for each $k$ shell.

The coefficients $S^{ij}_n(k_q)$ are evaluated by applying the
fast Fourier transform (FFT) to the template data along the $\psi$ grid
direction, requiring $O(\kmax^4 \log \kmax)$ work.

\section{Projection matching} \label{matchingsec}

Given an experimental image $\Mc \in \Mcb$ (or more precisely its 2D Fourier transform
$\Mchat(k,\psi)$ discretized on a polar grid), we seek to rank the templates
defined by $(\alpha_i,\beta_j)$ and the rotational degree of freedom
$\gamma$ in terms of how well they match the image.
For this we need a generative model for images.
We will then present an algorithm for matching $\gamma$.

\subsection{Image model with CTF correction}

For reasons having to do with the physics of data acquisition, each particle 
image $\Mc^{(m)}$ is not simply a projection of the electron-scattering
density, but may be modeled as the projection of the 
electron-scattering density convolved with a contrast transfer function
(CTF) \citep{Mindell2003}, plus noise.
The CTF includes diffraction effects due to the particle's depth in
the ice sheet, linear elastic and inelastic scattering, and
detection effects.
In the simplest
case, the CTF is radially symmetric and its Fourier transform
is real-valued and of the form $C^{(m)}(k)$, depending only on $k$.
For the purposes of the present paper, we will assume that 
$C^{(m)}(k)$ is known for each experimental image. By the convolution theorem,
the CTF acts multiplicatively on the Fourier transform image.
Thus, the expectation for the $m$th Fourier image
given by beam direction $(\alpha,\beta)$ and in-plane rotation $\gamma$,
deriving from a $\kb$-space scattering density $F$, is
$$
\Mchat^{(m)}(k,\psi-\gamma) \;\approx\;  C^{(m)}(k) S_{\alpha,\beta}[F](k,\psi)
~.
$$
To this is added measurement noise, which is assumed to be Gaussian
and i.i.d.\ on each pixel in the image.

\subsection{Fitting the best orientation for each image}
\label{s:match}

It would follow from the above noise model that the correct quantity to
minimize when searching for the best match would be the $L_2$ norm
of the difference between the template and the image,
which is equivalent to the $L_2$ norm
of their difference in the Fourier image plane.
This may be interpreted as minimizing the negative log likelihood,
i.e. finding the maximum likelihood.
Let $\Mchat = \Mchat^{(m)}$ be the current experimental image
in question, and $C = C^{(m)}$ be its corresponding known CTF.
We will denote by $\Mchat_\gamma$ the rotated Fourier image
\[
\Mchat_\gamma(k,\psi) :=
\Mchat(k,\psi-\gamma)~.
\]
{For the template index pair $(i,j)$,
a 
quantity to be minimized over the rotation $\gamma$ would thus be}
\begin{equation}
\| \Mchat_\gamma - 
C S_{\alpha_i,\beta_j}[F] \|_2,
\label{templatematch}
\end{equation}
where, for a function $g(k,\gamma)$ in polar coordinates, the $L_2$ norm
(up to a band-limit of $K$) is defined by
\be
\|g\|^2 := \int_0^K \int_0^{2\pi} |g(k,\psi)|^2 d\psi k dk.
\label{nrm}
\ee

However, there is typically uncertainty about the
overall normalization (a multiplicative prefactor)
associated with experimental images (as discussed for instance by Scheres
et al.\ \cite{Scheres2009}).
It is straightforward to check that
minimizing \eqref{templatematch} over $i,j,\gamma$ including an unknown
normalization factor for the image is equivalent
to maximizing over $i,j,\gamma$ the normalized inner product
\be
\frac{ \langle C S_{\alpha_i,\beta_j}[F], \Mchat_{\gamma} \rangle }
{ \|C S_{\alpha_i,\beta_j}[F]\| \, \| \Mchat_{\gamma} \| } ~.
\label{ip}
\ee
Here the inner product corresponds to the norm \eqref{nrm}, and
\eqref{ip} may be interpreted as the cosine of the ``angle'' between the
image and the template in the abstract vector space with norm \eqref{nrm}.
Maximizing \eqref{ip} over $i,j$ and $\gamma$
we refer to as {\em projection matching}.

For each $k_q \in \{k_1,\dots,k_{N_r} \}$, we precompute a Fourier series
representation of the image,
with $2N_q +1$ Fourier modes on the $q$th ring:
\[ \Mchat(k_q,\psi) = \sum_{n=-N_q}^{N_q} \Mchat_n(k_q) e^{i n \psi} . \]
This requires $O(\kmax^2 \log \kmax)$ work, done once
for each experimental image.

To rank the template matches we loop over all projection directions indexed
by $i,j$.
For each of these directions $(\alpha_i,\beta_j)$, since the denominator
of \eqref{ip} is fixed, 
we need to maximize the inner product
as a function of $\gamma$,
\[ 
\langle C S_{\alpha_i,\beta_j}[F], \Mchat_{\gamma} \rangle  :=
\int_0^K \int_0^{2 \pi}  \Mchat(k,\psi - \gamma)
\overline{C(k)\,S_{\alpha_i,\beta_j}[F](k,\psi)} d\psi k dk .
\]
Elementary Fourier analysis shows that
\[
\langle C S_{\alpha_i,\beta_j}[F], \Mchat_{\gamma} \rangle  =
\sum_{n} c_n(K) e^{-i n \gamma}
\]
where 
\be
c_n(\kappa) := 
\int_0^\kappa   \Mchat_n(k) \overline{C(k) S^{ij}_{n}(k)} k dk .
\label{cnk}
\ee
This integral is approximated using the existing $k$ grid
up to the maximum frequency $\kappa =O(\kmax)$, requiring $O(\kmax^2)$
effort to compute the $2N_q+1$ coefficients.
Since there are $O(\kmax^2)$ pairs $i,j$, this totals $O(\kmax^4)$ work per
experimental image.
Finally, for each pair $i,j$, we will compute the best value for $\gamma$ 
by tabulating $\sum_{n} c_n(K) e^{-i n \gamma}$ on a uniform grid of
$N_\phi$ values, using the FFT,
and taking the index with the maximum modulus. This requires
$O(\kmax^3 \log \kmax)$ work per image.
Thus, the complexity for angle fitting all images is $O(M\kmax^4)$.

\begin{definition}
The best values of $(\alpha,\beta,\gamma)$ for image $\Mc^{(m)}$
will be denoted by $\Ac_m = (\alpha_m,\beta_m,\gamma_m)$.
\end{definition}

In practice we implement two useful adjustments to the above procedure:

1) We search first over a coarse grid of template directions (with
angle resolution 5 times coarser than the grid defined by $N_\theta$
and $N_\phi$), and then only search over the set of grid points that
are within one coarse grid point of the global maxima found on the
coarse grid. This improves efficiency by a constant factor, and in our
tests does not degrade accuracy noticeably.

2) We define a ``randomization parameter'' $f_{rand}$.
If $f_{rand}=0$, we return the single best-fitting orientation:
$(\alpha_m,\beta_m,\gamma_m)$ defined above.
However, if $f_{rand}>0$, 
we instead choose $\Ac_m$ randomly from the set of all discrete
orientations $(\alpha,\beta,\gamma)$
that produce a normalized
inner product \eqref{ip} greater than $1-f_{rand}$.
Typically we choose $f_{rand}$ small, e.g. 0.02.
This uniformizes the distribution of orientations over the set which
fit the image almost equally well.
However, it is sometimes beneficial to increase $f_{rand}$ for improved convergence rate in the
least-squares procedure of the following section;
see section~\ref{resultsec2}.

It should be noted that
the idea of using Fourier methods to find the optimal third Euler
angle $\gamma$ can be found in  
\citep{Cong2003}, and the complexity of various alignment schemes
using polar grids in physical space is discussed in 
\citep{Joyeux2002}.

\begin{rmk}
In the simplest version of the frequency marching scheme we present,
the $c_n$ coefficients in \eqref{cnk} are evaluated using all
frequencies from zero to $\kappa$, since $F$, and hence the templates,
get updated over this frequency range each iteration.
We believe that updating the model $F$ only in the current $k$ 
shell could result in a faster algorithm without sacrificing accuracy.
In that case, the $c_n$ values are already known for a smaller $\kappa$,
and they can be incremented according to the formula
\[ c_n(\kappa_2) = c_n(\kappa_1) + 
\int_{\kappa_1}^{\kappa_2} \Mchat_n(k) \overline{C(k) S^{ij}_{n}(k)} k dk \, .
\]
\label{r:shell}
\end{rmk}

\section{Reconstruction from particle images with known angular assignments} 
\label{lsqsec}

Suppose now that we have the
collection of Fourier transforms of all $M$ particle images,
which we denote by
\[
\Mcbh \; :=\; 
\{ \Mchat^{(1)}(k,\psi),
\Mchat^{(2)}(k,\psi),
\dots,
\Mchat^{(M)}(k,\psi) \},
\]
together
with the corresponding angular assignments
$\Ac = \{\Ac_1,\Ac_2,\dots,\Ac_M \}$, where 
$\Ac_m = (\alpha_m,\beta_m,\gamma_m)$.

\begin{rmk}
In our present implementation, the values of $(\alpha_m,\beta_m,\gamma_m)$
are selected from the set of computed template angles, as discussed in the previous
section. A more sophisticated
projection matching procedure 
could generate ``off-grid" values for $\alpha, \beta$, and $\gamma$. The
reconstruction procedure below works in either case.
\end{rmk}

We seek to reconstruct a (new) model $F(k,\theta,\phi)$ that is consistent
with all of the data in a least squares sense. For each fixed spherical shell
of radius $k$, we use the representation \eqref{shellrep}, and
data $\Mchat^{(m)}_{\gamma_m}(k,\psi)$ which, if the angle assignments
are correct, is expected to approximate the CTF-corrected central slice
$C^{(m)}(k) S_{\alpha_m,\beta_m}[F](k,\psi)$ for the new model.
The matching between the new model and the image data on the $k$ shell
is done at a set of $N_\phi$ discrete angles $\phi_l$ on each great circle.
The Cartesian coordinates of those points are 
\begin{equation}
\left( \begin{array}{c} 
k [\cos \beta_m \cos \alpha_m \cos \psi_l - \sin \beta_m \sin \psi_l] \\ 
k [\sin \beta_m \cos \alpha_m \cos \psi_l + \cos \beta_m \sin \psi_l] \\ 
k [- \sin \alpha_m \cos \psi_l] 
\end{array} \right),
\label{datapts}
\end{equation}
from which their spherical coordinates
$(\theta_{ml},\phi_{ml})$ 
are easily computed.
Thus, for the Fourier sphere of radius $k$,
each image contributes $N_\phi$ points, for a total of $N_{tot} = M N_\phi$ points
on this sphere.

At this stage, we simply collect all image data
with given $k$
into a right-hand side vector $\db$, defined by its $N_{tot}$ entries
$$
d_{(m-1)N_\phi + l} = \Mchat^{(m)}_{\gamma_m}(k,\psi_l)
~,\qquad \mbox{ for } m=1,\dots,M, \quad l=1,\ldots,N_\phi ~.
$$
From now we shall use a single index $i = 1,\ldots,N_{tot}$ to reference
these entries, and use $(\theta_{i},\phi_{i})$ for the spherical
coordinates for the corresponding point on the great circle.
(In our present implementation
we ignore the use of image normalization factors
that could be extracted from the orientation fitting algorithm presented in
section~\ref{s:match}. This is justified since our synthetic
noisy images will be generated without varying normalization factors.)

Let $\fb$ denote the vector of ``unrolled" spherical harmonic
coefficients in the representation for $F$ on the current $k$ sphere,
\[ \fb = \{f_{0,0}(k),f_{1,-1}(k),
f_{1,0}(k),f_{1,1}(k),f_{2,-2}(k),\dots,f_{p(k),p(k)}(k) \}~,
\]
so that $\fb$ is a complex vector of length $(p(k)+1)^2$.
Then we define the complex-valued
matrix S, of dimension $N_{tot} \times (p(k)+1)^2$, by 
\begin{equation}
\label{Sdef}
(S \, \fb)_i  = 
 \sum_{n=0}^{p(k)} \sum_{m=-n}^n f_{nm}(k) Y_n^m(\theta_i,\phi_i) \, .
\end{equation}
Thus, $S$ evaluates the spherical harmonic expansion with coefficients
$\fb$ at all sphere points.
We also let $C_i(k)$ be the corresponding CTF
for each image at frequency $k$ (note that $C_i(k)$ is the same
for all $i$ belonging to the same image), and let $\Cc$ be
the diagonal matrix with diagonal entries $C_i(k)$.
We find the desired solution by solving the problem 
\[ \Cc S\fb  = \db
\] 
in a least squares sense.
Note that the use of the $l^2$ norm
corresponds to the same assumption of i.i.d.\ Gaussian noise
on the images as in the previous section.

For this, we use conjugate gradient (CG) iteration on the normal equations
\[ 
S^H \Cc^H \Cc  S \, \fb = S^H \Cc^H \db,
\] 
where $S^H$ is the Hermitian transpose of $S$.
Applying the diagonal matrices $\Cc$ and $\Cc^H$ is trivial.
Direct computation of the matrix-vector product with $S$ or $S^H$ would require
$O(M \, N_\phi  \, (p(k)+1)^2) = O(M \kmax^3)$ work.
The application of $S$ is easily accelerated, however, by first evaluating
the spherical harmonic expansion induced by $\fb$ on a {\em regular} grid with
$N_\phi$ equispaced points in the $\phi$ direction, and
$N_\theta$ points in the $\theta$ direction whose cosines are 
classical Chebyshev nodes (see section \ref{sec:disc}).
Using separation of variables, this requires
$O(\kmax^3)$ work, since $N_\phi = N_\theta = O(\kmax)$.
The value of $S\,\fb$ at an arbitrary point can then be computed by 
local $q$th-order
Lagrange interpolation in $\phi$ and $\theta$ from this
regular grid, at a cost of $q^2$ operations
per target. It is sufficient to use $q$ in the range 5 to 10.
We denote by $S_{reg}$ the mapping from spherical harmonic coefficients
to values on the $N_\phi \times N_\theta$ grid, and by $T$ the 
(sparse) interpolation matrix from the 
$N_\phi \times N_\theta$ grid points to the $N_{tot}$ arbitrary locations.
In other words, up to interpolation error,
\[ S \fb \approx T S_{reg} \fb, \]
and applying $S$ in this manner requires only $O(M q^2 \kmax + \kmax^3)$ work.
Clearly,
\[ 
S^H \Cc^H \Cc S \approx S^H_{reg}  T^H \Cc^H \Cc T S_{reg}, 
\]
where $S^H_{reg}$ can be applied by projection from a regular grid to
a spherical harmonic expansion, which is also $O(\kmax^3)$.
Thus each CG iteration requires
only $O(M q^2 \kmax + \kmax^3)$ work. So long as the $N_{tot}$ points
have
reasonable coverage of the sphere (specifically that there are
no large ``pockets'' on the sphere which are empty of points),
the system is well conditioned
and requires only a modest number of iterations, around 20--50,
to achieve several digits of accuracy.
It is worth noting (as in other Fourier-based 
schemes, such as FREALIGN \citep{Grigorieff2007,Grigorieff2016,Lyumkis2013}) 
that this least squares
procedure is the step in the reconstruction that is responsible for denoising. 
The more experimental images that
are available (with poor signal-to-noise ratios but correctly assigned angles), 
the more accurately we
are able to estimate $\fb$, hence $F(\kb)$ and the
electron-scattering density $f(\xb)$.

The above describes the computation
of $F$ on a single $k$ shell.
One 
advantage of our representation is that 
to build the complete new model $F$, the least squares solve for each
$k$ shell in the grid $\{k_1,\ldots,k_{N_r}\}$ may be performed independently.
This gives an overall
complexity for the reconstruction of $O(M q^2 \kmax^2 + \kmax^4)$.

\begin{definition}
Let $\Fb$ indicate the entire set of model coefficients
$\{ \fb(k_1),\ldots,\fb(k_{N_r})\}$.
We will refer to its least squares approximation (on all $k$ shells) by $\Fb^*$. We summarize this by
\be
\Fb^* =L(\Mcbh,{\cal A}):=\arg\min_{\Fb} 
\sum_{m=1}^{M} \|\Cc^{(m)} S_{\Ac_m}(\Fb)- \Mchat^{(m)}\|_2^2 \, ,
\label{Fstardef}
\ee
where $\Cc^{(m)}$ is the CTF for the $m$th image
and, by analogy with \eqref{Pabg},
$S_{\Ac_m}(\Fb)$ denotes the slice
\[ S_{\Ac_m}[F] = \widehat{P_{\alpha_m,\beta_m,\gamma_m}[f]}~. \]
{Here}, $F$ is determined from the coefficients $\Fb$ 
via the spherical harmonic representation \eqref{shellrep}.
\end{definition}

\section{The full inverse problem}
\label{RLsec}

As discussed in the previous section, 
if the angular assignments ${\cal A}=\{\Ac_1,..,\Ac_M \}$ of the
experimental images were given, $\Fb$ could be recovered by solving the 
least squares problem \eqref{Fstardef}.
Since the angles ${\cal A}$ are unknown, however, it is standard
to write the cryo-EM reconstruction problem in terms of the nonlinear and
nonconvex optimization task in the joint set of unknowns:
\begin{equation}
\label{fullopt}
\{\Fb^*,\Ac^*\}=\arg\min_{\{\Fb,\Ac\}} 
\sum_{m=1}^{M} \|\Cc^{(m)} S_{\Ac_m}(\Fb)- \Mchat^{(m)}\|_2^2 \, .
\end{equation}
As noted above, for simplicity of presentation of the scheme, we omit
image normalization prefactors discussed in section~\ref{s:match}.

One well-known approach to solving \eqref{fullopt} is to start with an initial
guess $\Fb^{(0)}$ for the representation $\Fb$,
then to iterate as follows:

\vspace{1.5ex}  
\framebox{\begin{minipage}{15cm} 
{\bf Classical iterative refinement}

\vspace{1ex}
Set $i=1$. While (convergence criterion has not been met),
\begin{enumerate}
\item Compute $\Ac^{(i)}$ from $\Fb^{(i-1)}$ by projection matching
(section \ref{matchingsec}).
\item Set $\Fb^{(i)} = L(\Mcbh,\Ac^{(i)})$
by least squares solution (section \ref{lsqsec}).
\item $i \leftarrow i+1$.
\end{enumerate}
Compute $F$ from the final $\Fb$
via \eqref{shellrep} then take the inverse Fourier transform
\eqref{FTinvdef} to recover $f$.
\end{minipage}}
\vspace{1.5ex}

This iteration can be viewed as
coordinate descent, alternating between fitting the best
angle assignment and fitting the best model density.
It will converge, but to a local minimum---not necessarily the correct
solution \citep{Cheng2015}. Various attempts to overcome this convergence failure
have been proposed, including annealing strategies and stochastic hill climbing
(see, for example, \citep{Elmlund2012,Cheng2015}), but robustness has remained an issue.
 
\subsection{Frequency marching (recursive linearization)}

As in the iteration above, we alternate between projection matching to obtain
estimates for the angles $\Ac$
and solving least squares problems to determine the best set of density
coefficients $\Fb$.
However, by continuously increasing the 
resolution---measured in terms of the maximal spatial frequency used in the
representation for $\Fb$---we bypass the difficulties associated with multiple
minima in existing attempts at iterative refinement.
In the language of optimization, this can be viewed as a 
{\em homotopy method} using the maximum spatial frequency
(resolution) as the homotopy parameter.
The basic intuition underlying our scheme is motivated 
by the success of recursive linearization in inverse acoustic scattering 
\citep{BaoLi2015,Borges2016,Chenrep1088,Chen}.

More precisely, let $\Mcbh([0,k])$ denote the set of Fourier transforms
of all experimental images restricted to the disk of radius $k$, and
let $\Fb([0,k])$ denote the density coefficients only up to frequency
$k$, i.e.\ using the shells for which $k_q\le k$.
The full objective function minimization restricted to maximum frequency $k$ is
\begin{equation}
\{\Fb^*([0,k]),{\cal A}^*\}=\arg\min_{\{\Fb([0,k]),{\cal A}\}} 
\sum_{m=1}^{M} \| \Cc^{(m)} S_{\Ac_m}(\Fb([0,k]))-\Mchat^{(m)}([0,k]) \|_2^2~,
\end{equation}
which is still a nonlinear and nonconvex optimization problem.
However, if $\Fb([0,k])$ is known and we only seek to find
$\Fb([0,k+\delta k])$ for sufficiently small $\delta k$, then the 
solution can be reached by a suitable linearization.
Moreover, at low frequency, say for $k \leq k_1 = 2$, the 
landscape is extremely 
smooth so that the global minimum is easily located. Thus, we propose simply
assigning random angles 
$\Ac_m = (\alpha_m,\beta_m,\gamma_m)$ to the images at $k_{min} = k_1 = 2$
and iterating as follows:

\vspace{1.5ex} 

\framebox{\begin{minipage}{15cm} 
{\bf Solution by frequency marching}

\vspace{1ex}
On input, we define a sequence of frequency steps
from $k_1$ to $k_{N_r} = \kmax$ with a step of $\delta k = k_{i+1} - k_i$.

$\bullet$ Set $\Ac^{(0)}$ to uniform random values over the allowed ranges. \\
$\bullet$ Set $\Fb([0,k_1]) = L(\Mcbh([0,k_1]),\Ac^{(0)})$
by least squares (section \ref{lsqsec}).

\vspace{1ex}
For $i=1,2,\ldots,N_r-1,$
\begin{enumerate}
\item Compute $\Ac^{(i)}$ from $\Fb([0,k_i])$ by projection matching
(section \ref{matchingsec}).
\item Set $\Fb([0,k_{i+1}]) = L(\Mcbh([0,k_{i+1}]),\Ac^{(i)})$
by least squares (section \ref{lsqsec}).
\end{enumerate}
Compute $F$ from $\Fb([0,\kmax])$
via \eqref{shellrep} then take the inverse Fourier transform
\eqref{FTinvdef} to recover $f$.
\end{minipage}}
\vspace{1.5ex}

Note that the procedures of sections \ref{matchingsec}
and \ref{lsqsec} restrict naturally to any frequency range $[0,k]$.
The key feature of refinement by recursive marching is that it is 
a deterministic procedure involving only linear solves
and angular assignments of images.
The overall complexity, combining those from sections \ref{matchingsec} and
\ref{lsqsec}, and summing over $k$, is $O(M\kmax^5 + \kmax^6 + 
q^2M\kmax^3)$,
where the first two terms come from angle matching and the last from least
squares solution.
Since in current cryo-EM applications, $\kmax \sim 10^2$, while
$M\sim 10^5$ to $10^6$, and $q^2<\kmax$, the dominant cost is
$O(M\kmax^5)$, 
assuming we use the global angle matching procedure described above.
However, we find that in our examples,
due to the number of CG iterations, the
time for least squares fitting is actually quite similar
to that for angle matching, i.e. they are close to balanced.

While we are not able to provide 
a proof in the general case, we believe that under fairly broad
conditions this iteration will converge with high probability 
for sufficiently small $\delta k$
(the radial frequency grid spacing).
Informally speaking, we believe that a solution near the 
global minimum is {\em often} reached at low $k$ in the first few 
iterations and that a path to the global minimum at $\kmax$ is then 
reached as $k$ increases continuously.

In the experiments below we show that, even with very noisy data,
a harsh random start at
$k_{min} = 2$ is sufficient. We will return to this question in the 
concluding section.

In practice, we have implemented a further acceleration to the above scheme:
if the mean absolute change in angles between $\Ac^{(i+1)}$ and $\Ac^{(i)}$
is less than $10^{-3}$,
the next increment of the index $i$ is set to five rather than its usual
value of one.
This greatly reduces the number of steps in the marching procedure
once the majority
of angles have locked in with sufficient accuracy.

\section{Numerical experiments with synthetic data}
\label{resultsec1}

\begin{figure}
\begin{tabular}{ccc}
Rubisco
&
Lipoxygenase-1
&
Neurotoxin
\\
\includegraphics[width=.28\linewidth,trim= {0cm 4cm 1cm 0cm},clip]{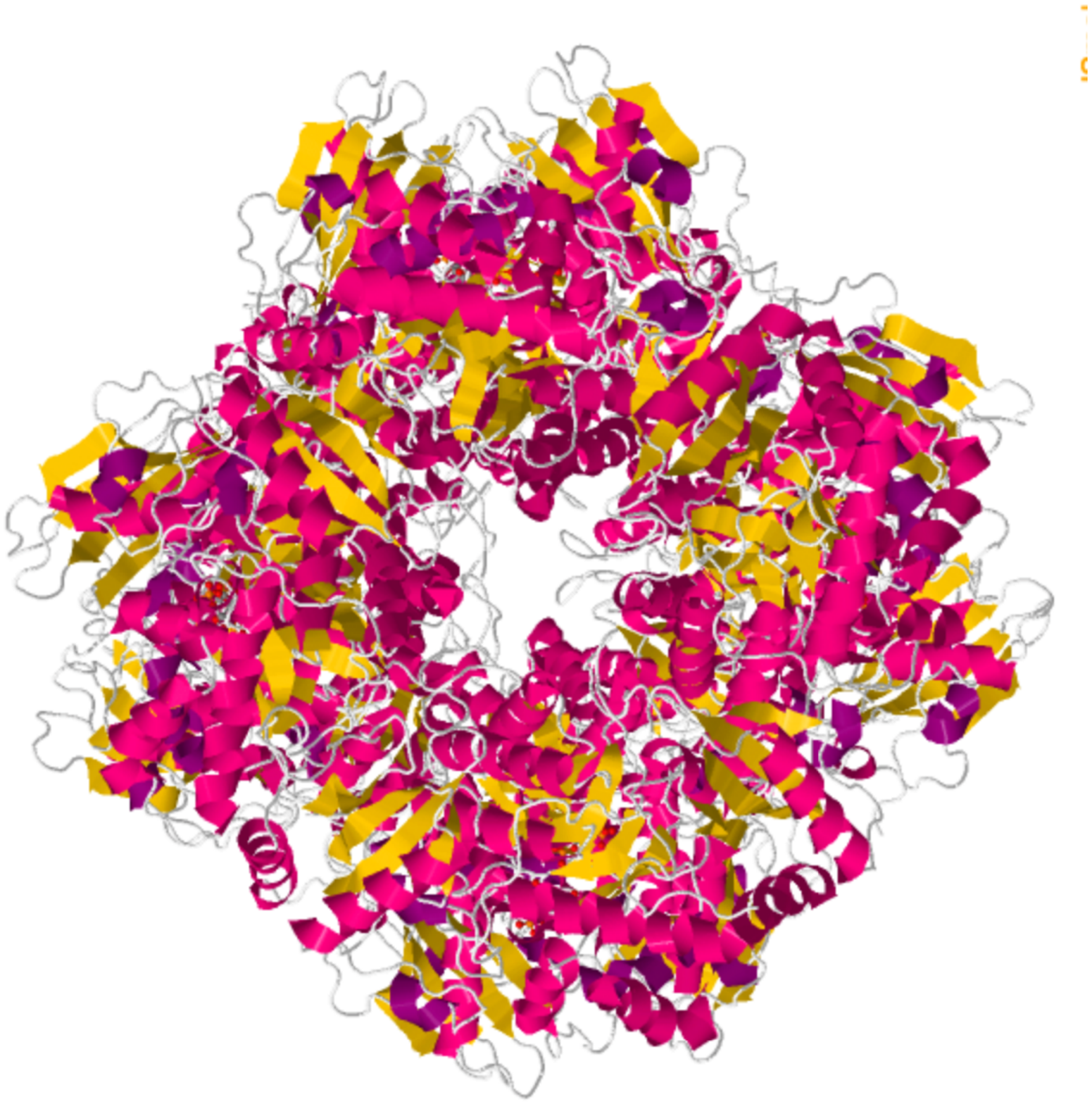}
&
\includegraphics[width=.33\linewidth,trim= {0cm 4cm 1cm 0cm},clip]{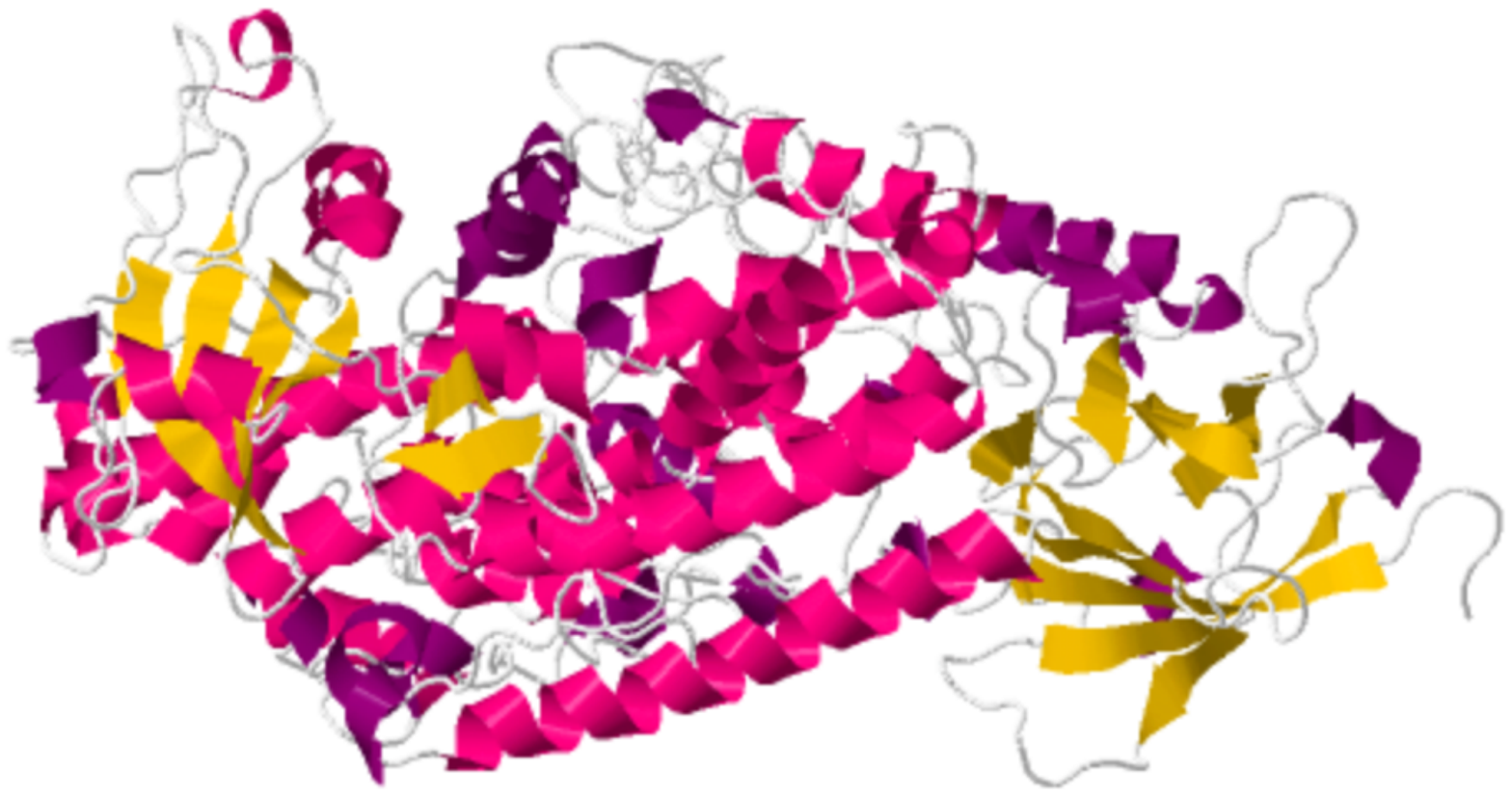}
&
\includegraphics[width=.33\linewidth,trim= {0cm 4cm 1cm 0cm},clip]{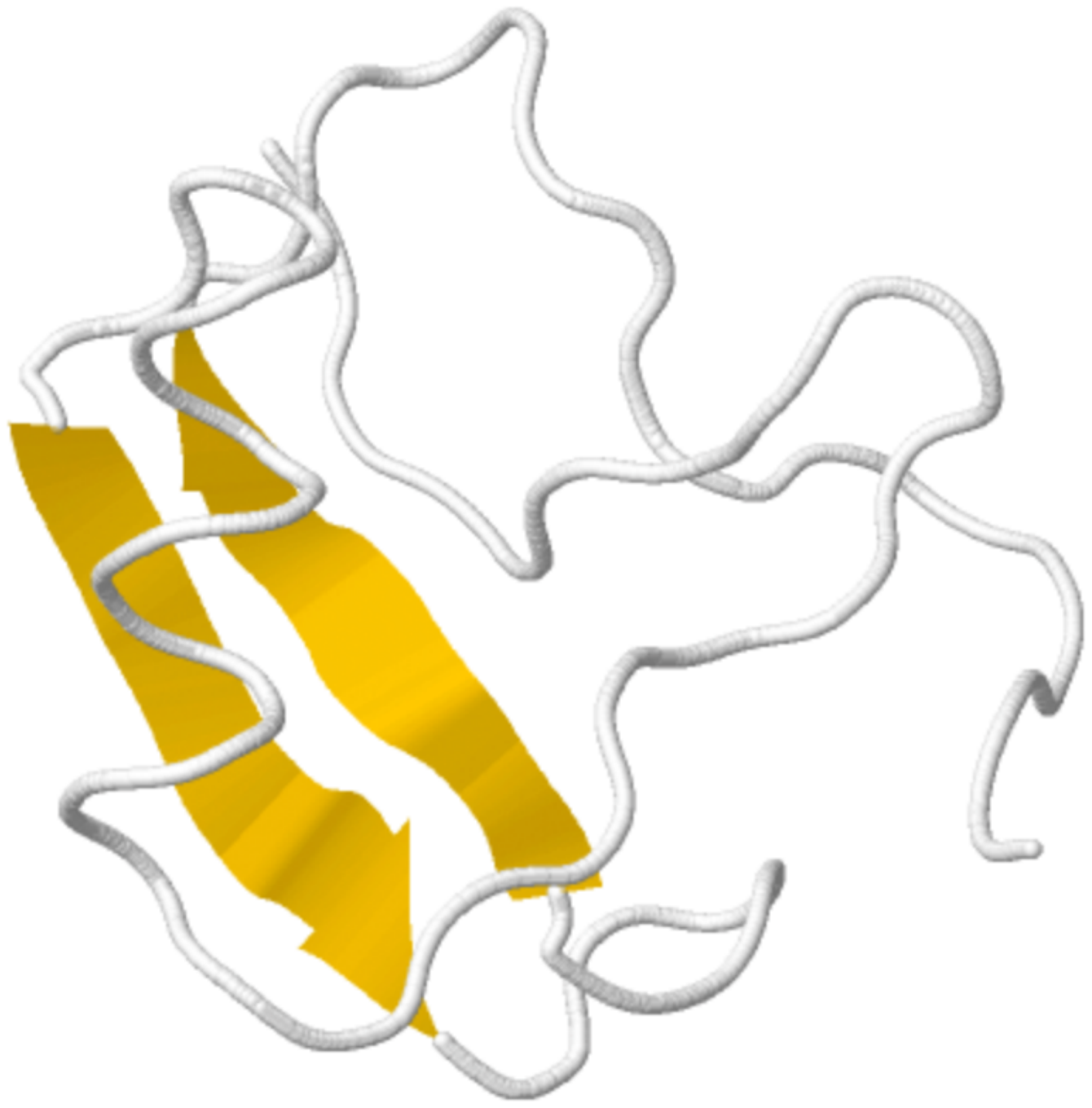}
\\
\includegraphics[width=.28\linewidth,trim= {6cm 9cm 6cm 9cm},clip]{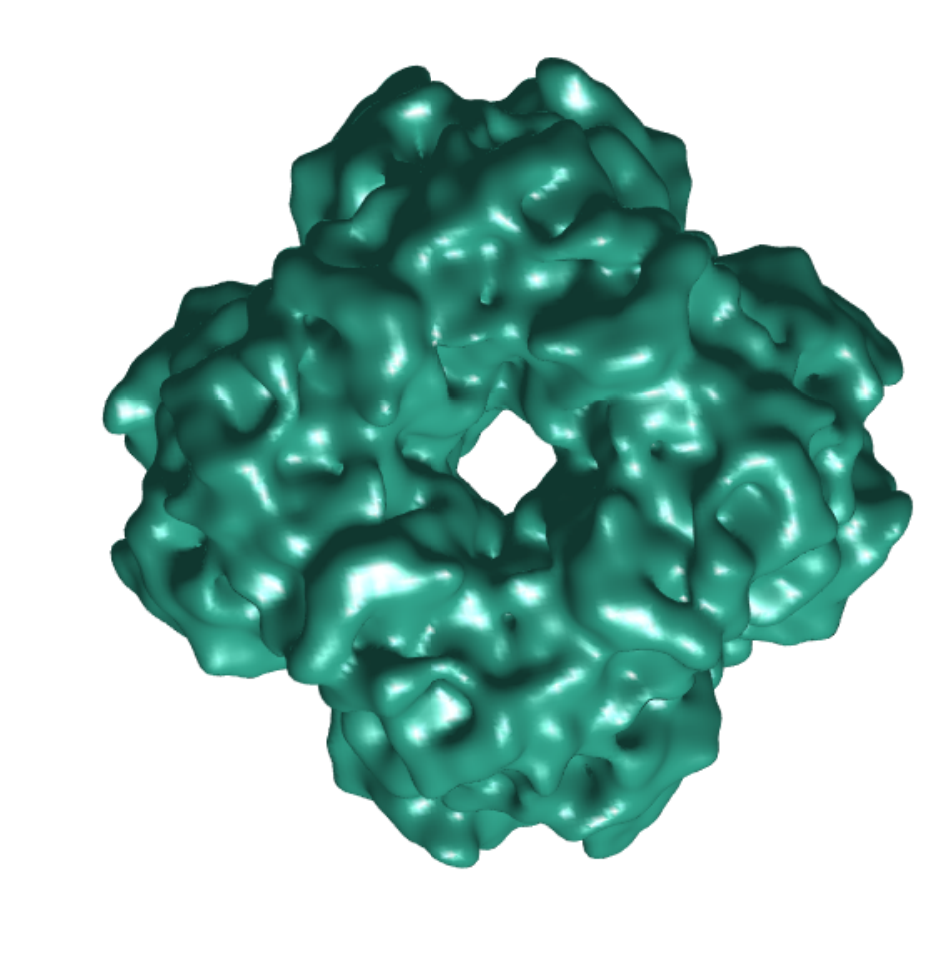}
&
\includegraphics[width=.33\linewidth,trim= {5cm 10cm 5cm 10cm},clip]{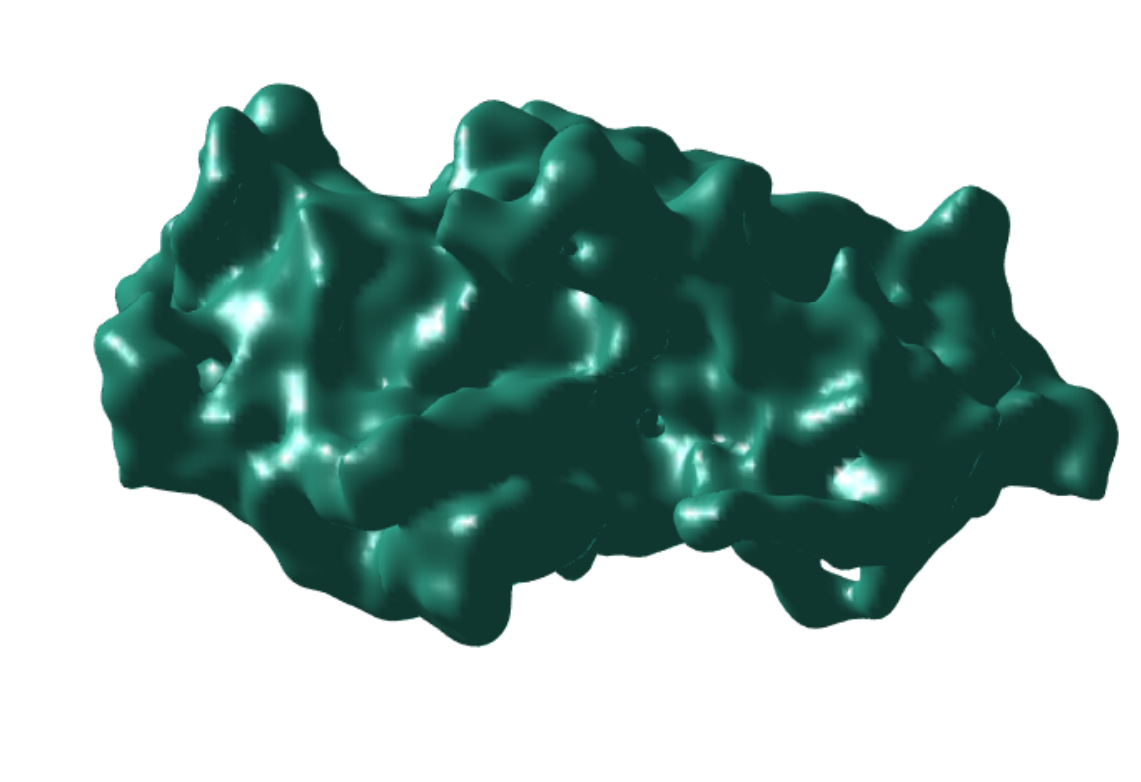}
&
\includegraphics[width=.33\linewidth,trim= {5cm 8cm 5cm 8cm},clip]{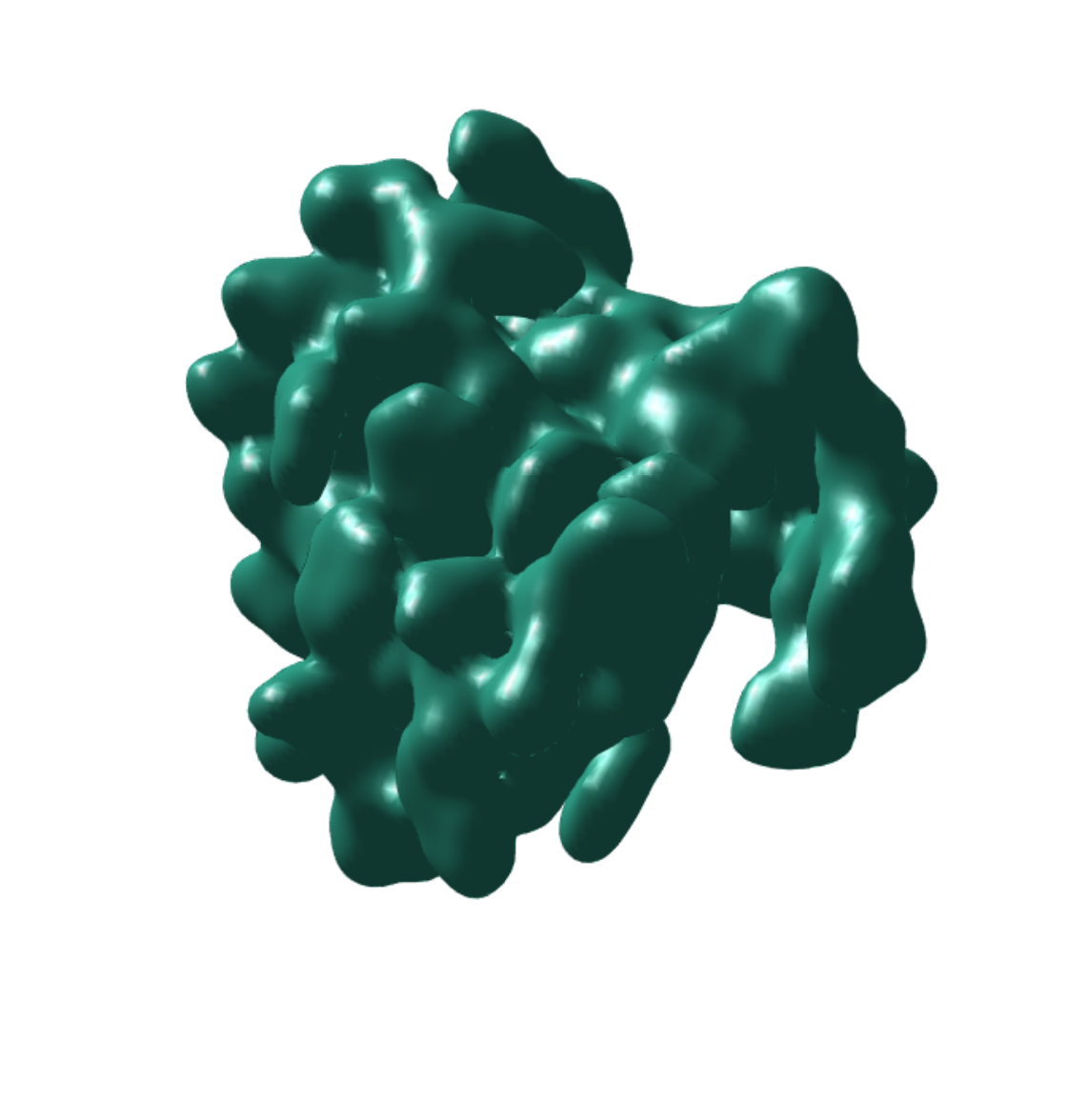}
\\
\end{tabular}
\caption{Top row: The three protein
structures (derived from X-ray crystallography). Bottom row: the
  corresponding atomic densities $f(\xb)$ used as a ground truth for the
  numerical experiments,
shown via an isosurface at approximately half the maximum density.
Note that the visual length scale for the three molecules is
different.
}
\label{figure:three_models}
\end{figure}

We implemented the above frequency marching algorithm
and performed experiments using simulated images to reconstruct
electron scattering densities for three proteins.

\subsection{Choice of ground-truth densities and error metric}

The three proteins we used for numerical experiments were (see
Fig.~\ref{figure:three_models}): spinach rubisco enzyme (molecular
weight 541kDa, 133 \AA\ longest dimension), lipoxygenase-1 (weight
95kDa, 107 \AA\ longest dimension), and scorpion protein neurotoxin
(weight 7.2kDa, 33 \AA\ longest dimension).  Atomic locations in
\AA\ are taken from the protein data bank (PDB) entries \cite{pdb},
with codes 1RCX, 1YGE, 1AHO, respectively, but were shifted to have
their centroid at the origin.  Only in the last of the three are
hydrogen atom locations included in the PDB. Thus, the number of atoms
used were 37456, 6666, and 925, respectively.  Note that the last of
these, the neurotoxin, is much smaller than the $\approx 100$ kDa
lower weight limit that can be imaged through cryo-EM with current
detector resolution and motion-correction technology; our point in
including it is to show that, given appropriate experimental image
resolution, reconstruction of small molecules does not pose an
algorithmic problem.

For each protein, a ground-truth
electron scattering density $f(\xb)$ was produced using the following
simple model.
Firstly, a non-dimensionalization
length $D$ was chosen giving the number of physical \AA\ per unit
numerical length, such that the support of $f$ in such numerical
units lies within the unit ball.
For rubisco, $D = 70$ \AA, for lipoxygenase-1, $D = 60$ \AA, and
for the neurotoxin, $D = 25$ \AA.
Thus, all physical distances in what follows are divided by $D$
in the numerical implementation.
We summed 3D spherically-symmetric
Gaussian functions, one at each atomic location,
giving each Gaussian a standard deviation $\sqrt{(0.5 r)^2 + b^2}$,
where $r$ is the atomic
radius for the type of atom (which vary from 0.42 \AA\ to 0.88 \AA).
The factor $0.5$ results in around 74\% of the mass of the
unblurred Gaussian falling within radius $r$,
and $b$ is a convolution (blurring) radius used to make $f$ smooth enough
to be accurately reconstructed using the image pixel resolution.
The radii $b$ used were 2.5\AA\ for rubisco,
2.0\AA\ for lipoxygenase-1, and 1.0\AA\ for scorpion toxin,
chosen to be around twice the simulated pixel spacing (see next section).
Thus, when we assess reconstruction errors, we are doing so against
an appropriately smoothed ground-truth $f$.
For simplicity, each such Gaussian was given a unit peak amplitude.

Our error metric for a reconstructed density
$\tilde f$ is the relative $L_2$-norm,
\be
\epsilon := \frac{\|\tilde f - f\|_{L_2([-1,1]^3)}}{\|f\|_{L_2([-1,1]^3)}}
~,
\label{err}
\ee
which is estimated using quadrature on
a sufficiently fine uniform 3D Cartesian grid (we used $100$ points in
each dimension).

\begin{rmk}
Since $\tilde f$ generally acquires an arbitrary rotation
$(\alpha,\beta,\gamma)$ relative to the ground-truth $f$, we must 
first rotate it to best fit the original $f$.
This is done
by applying the procedure of section~\ref{matchingsec}
to a small number (typically 10) of random projections.
$\tilde f$ is then rotated using the 3D non-uniform FFT \cite{nufft}
to evaluate its Fourier transform $\tilde F$ at a
rotated set of spherical discretization points as in section~\ref{sec:disc}.
A second non-uniform FFT  is then applied 
to transform back to real space.
Finally the error \eqref{err} is evaluated.
\end{rmk}

Some researchers use a normalized cross-correlation metric to report
errors (eg.\ \cite{Joubert2015}). We note that, when errors are
small, this metric is close to $1 - \epsilon^2/2$ with $\epsilon$
given by \eqref{err}.

\subsection{Generation of synthetic experimental images}

For each of the three proteins, $M$ (of order 50,000)
synthetic experimental images were produced as follows.
For each image $m$ we
first defined a realistic radially-symmetric CTF function $C^{(m)}(k)$
using the standard formulae \cite{Mindell2003,wade}
\be
C(k) = B(\theta) [w_1 \sin\chi(\theta) - w_2 \cos \chi(\theta)],
\quad \mbox{ where } \;
\chi(\theta) := \frac{1}{2} k z \theta^2 + k C_s \theta^4 / 8,
\quad
B(\theta) := e^{-\theta^2/\theta_0^2}~.
\ee
Here $\chi$ is called the phase function, $\theta_0 = 0.002$
sets the microscope acceptance angle, $w_2=0.07$ controls the
relative inelastic scattering, with $w_1^2+w_2^2=1$,
and the spherical aberration is $C_s = 2\times 10^7$ \AA.
The defocus parameter $z$ was different for each image,
chosen uniformly at random in the interval $[1,4] \times 10^4$ \AA.
Finally, in the above,
angles $\theta$ are related to wavenumbers $k$
in the numerical experiments via
\be
\theta = \frac{\lambda k}{2\pi D} ~,
\ee
where $\lambda = 0.025$ \AA\ is the free-space electron wavelength
(a typical value corresponding to a 200 keV microscope),
and for convenience the distance scaling $D$ has been included.
In our setting $\theta$ is of order $10^{-5}$ times the numerical wavenumber $k$.

The images were sampled on a uniform 2D grid at a standard
resolution of $100\times 100$ pixels
covering the numerical box $[-1,1]^2$.
The pixel spacing thus corresponded to $D/50$, or between 0.5 \AA\ for
neurotoxin (this is smaller than currently achievable experimentally)
and 1.4 \AA\ for rubisco.
The noise-free signal images were produced by inverse Fourier transforming
(via the 2D non-uniform FFT) slices taken at random
orientations through the ground-truth Fourier density $F(\kb)$,
after multiplication by the radial CTF $C^{(m)}(k)$ particular to each image.
Finally, i.i.d.\ Gaussian noise was added to each pixel,
with variance chosen to achieve a desired signal-to-noise ratio (SNR).
SNR has the standard definition in imaging as the ratio of the
squared $L_2$-norm of the signal to that of the noise
(for this we used the domain $[-1,1]^2$ since the molecule images
occupy a large fraction of this area).
We generated images at SNR values of
$\infty$ (no noise), 0.5, 0.1 and 0.05.
A typical value in applications is 0.1.

\subsection{Results}
\label{resultsec2}

All numerical experiments were run on desktop workstations with 14
cores, Intel Xeon 2.6GHz CPU, and 128 GB RAM.
Our implementation is in Fortran, with OpenMP parallelization.
We used frequency marching
with step size $\delta k = 2$, starting 
{with random angular assignment
at $k_{1} = 2$, and marching to full resolution at
$\kmax = 70$. This $\kmax$ was sufficient to resolve the decaying tails of $F$
without significant truncation error, given our choice of blurring radius $b$ at
approximately two pixels.  We used interpolation order $q=7$, and we
started with randomization factor $f_{rand} = 0.02$.
The latter was changed adaptively during marching to balance solution time and accuracy, as follows:
if at some point, the least squares CG did not converge after 100 iterations,
$f_{rand}$ was doubled and the least squares solve repeated.
If instead CG required less than 50 iterations, $f_{rand}$ was halved.
The speed of our algorithm allowed us to carry out 
multiple runs for each protein, using either the same or a
fresh set of synthetic images.}
{In the present paper, we carried out 5 runs for each protein.}

To assess how close our relative
$L_2$ error $\epsilon$ was to the best achievable, given
the image sampling, Fourier representation, and SNR, we
computed for comparison
the {\em best possible} $\epsilon$
achieved by reconstructing $f$ through a single application
of {the least squares procedure} (section \ref{lsqsec})
with all image orientations set to their true values.
Table \ref{table-errors} shows that the errors resulting from applying our
proposed frequency marching algorithm
{exceed this best possible error by only around $10^{-2}$, for all molecules and noise levels.}

\begin{figure}
\bc
\begin{tabular}{cc}
SNR 0.5 
& 
SNR 0.05 
\\ 
(a)\raisebox{-1in}{
\includegraphics[width=.25\linewidth,trim= {5cm 9cm 5cm 9cm},clip]{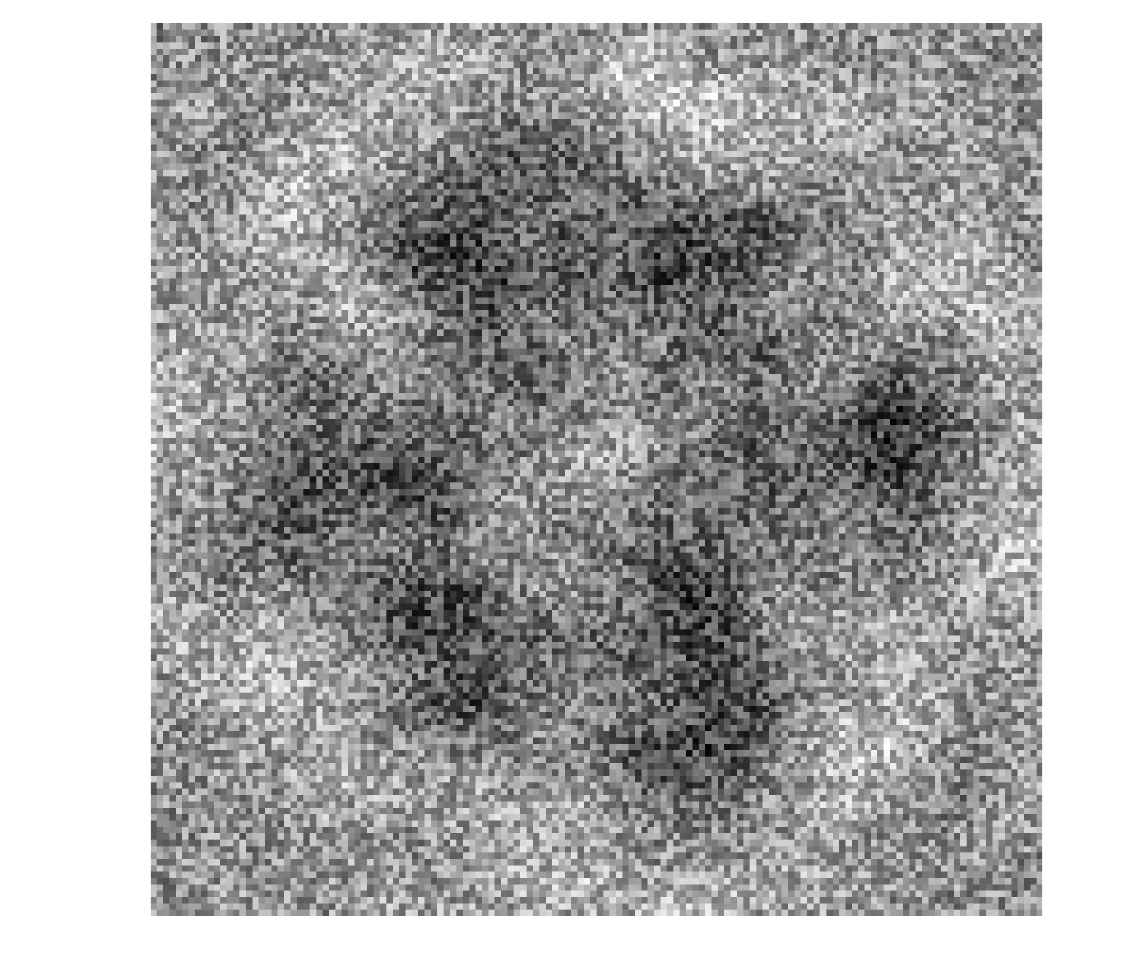}}
&
(b)\raisebox{-1in}{
\includegraphics[width=.25\linewidth,trim= {5cm 9cm 5cm 9cm},clip]{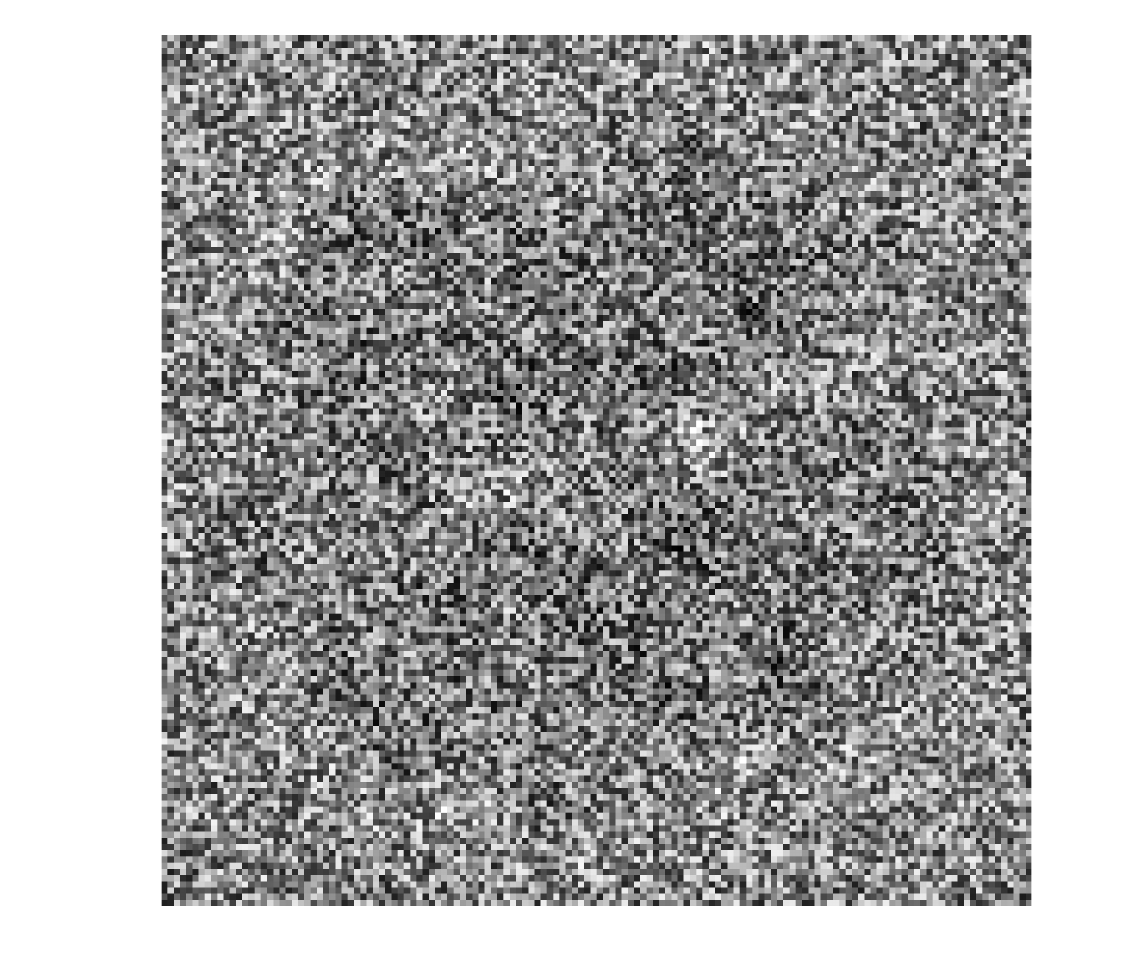}}
\\
(c)\raisebox{-1in}{
\includegraphics[width=.25\linewidth,trim= {7cm 11cm 7cm 10.5cm},clip]{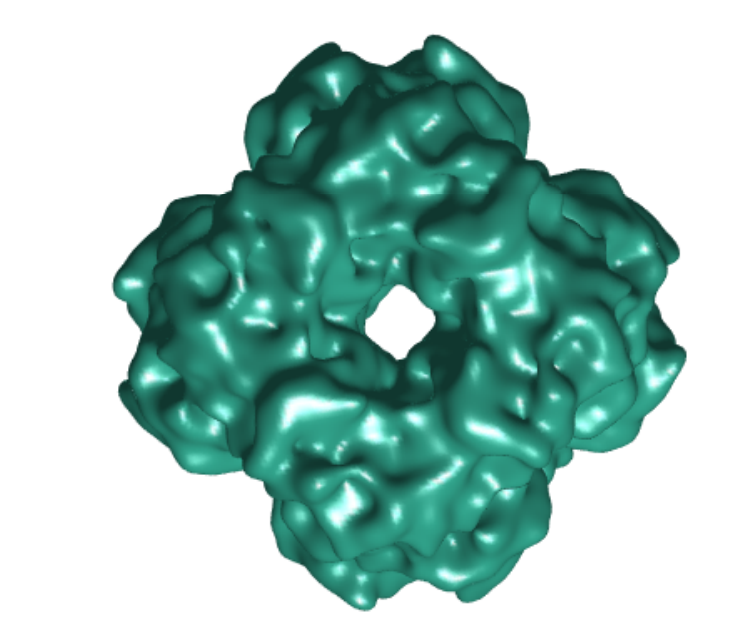}}
&
(d)\raisebox{-1in}{
\includegraphics[width=.25\linewidth,trim= {7cm 11cm 7cm 10.5cm},clip]{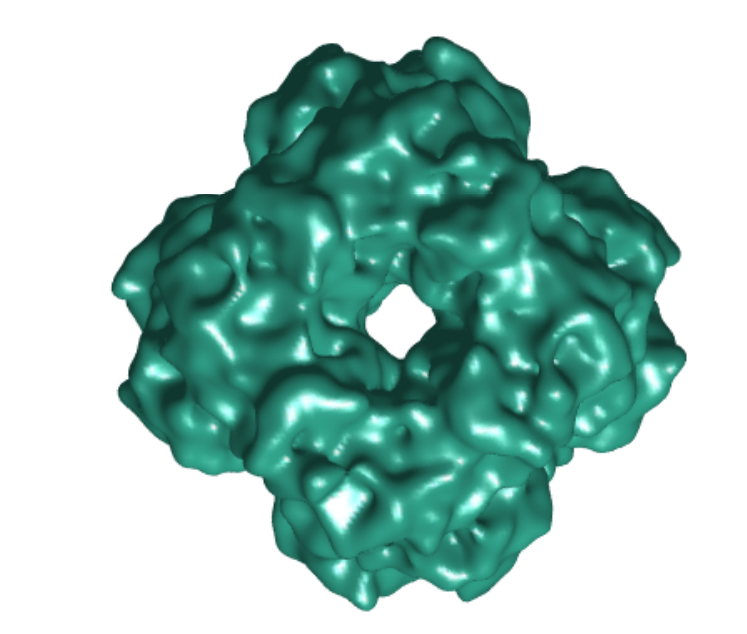}}
\\ 
(e)\raisebox{-1in}{
\includegraphics[width=.25\linewidth,trim= {5cm 9cm 5cm 9cm},clip]{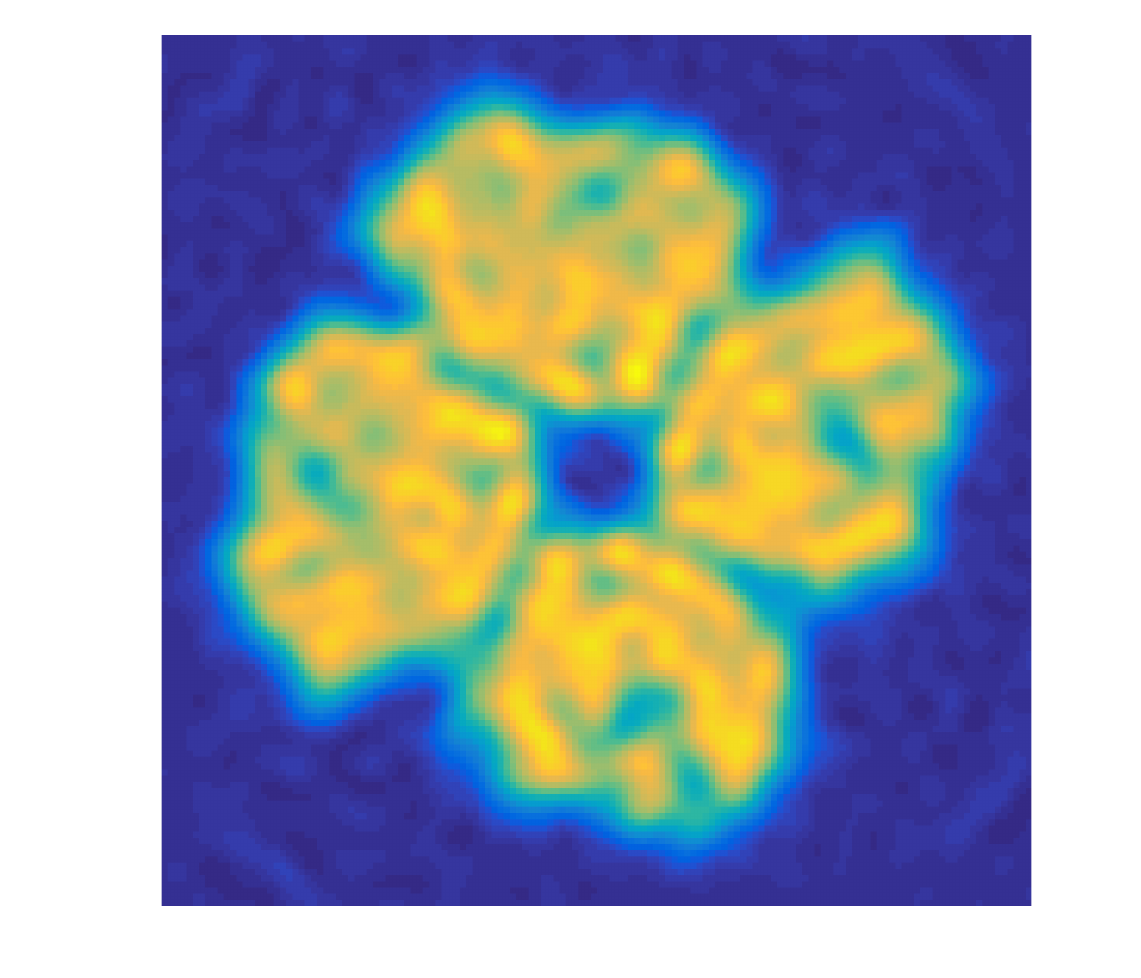}}
&
(f)\raisebox{-1in}{
\includegraphics[width=.25\linewidth,trim= {5cm 9cm 5cm 9cm},clip]{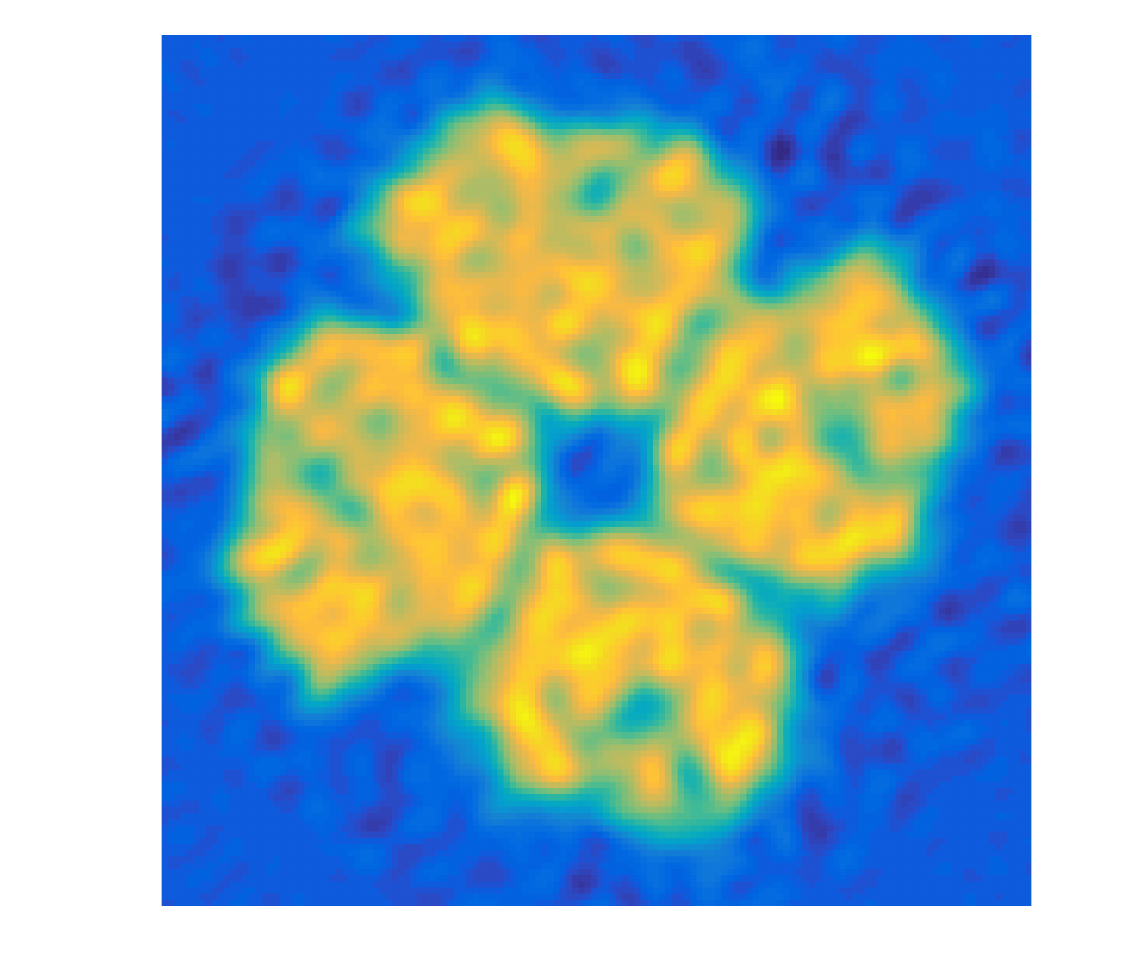}}
\end{tabular}
\ec
\caption{Results for Rubisco reconstruction. (a,b): Examples of
  experimental images at SNR 0.5 and 0.05 (black is smallest values, white largest). (c,d): Reconstructions
  by frequency marching using experimental images at SNR 0.5 and 0.05, shown as isosurfaces. (e,f): 
Examples of a single slice through the reconstructed $f$ (blue is smallest values, yellow largest).}
\label{figure:rubisco-results}
\end{figure}

\begin{figure}
\bc
\begin{tabular}{cc}
SNR 0.5 
&
SNR 0.05 
\\
(a)\raisebox{-1in}{
\includegraphics[width=.25\linewidth,trim= {5cm 10cm 5cm 9cm},clip]{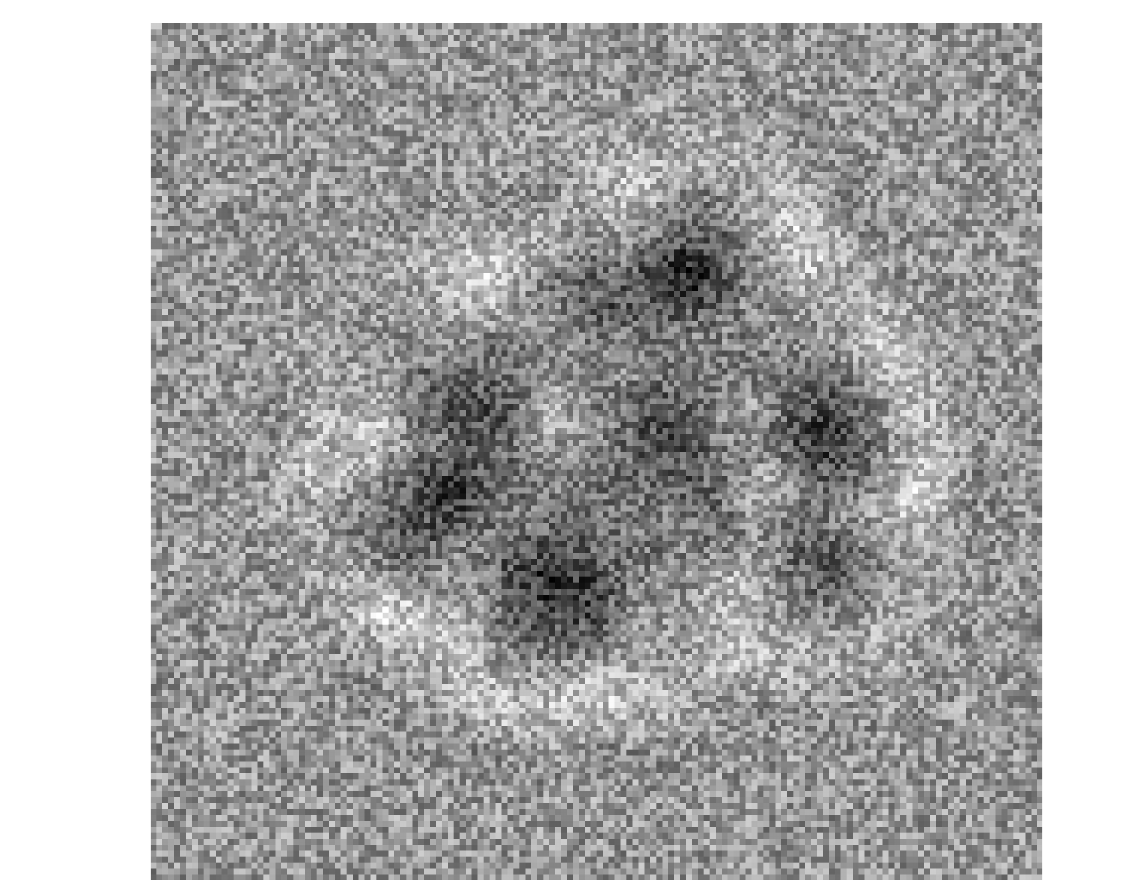}}
&
(b)\raisebox{-1in}{
\includegraphics[width=.25\linewidth,trim= {5cm 10cm 5cm 9cm},clip]{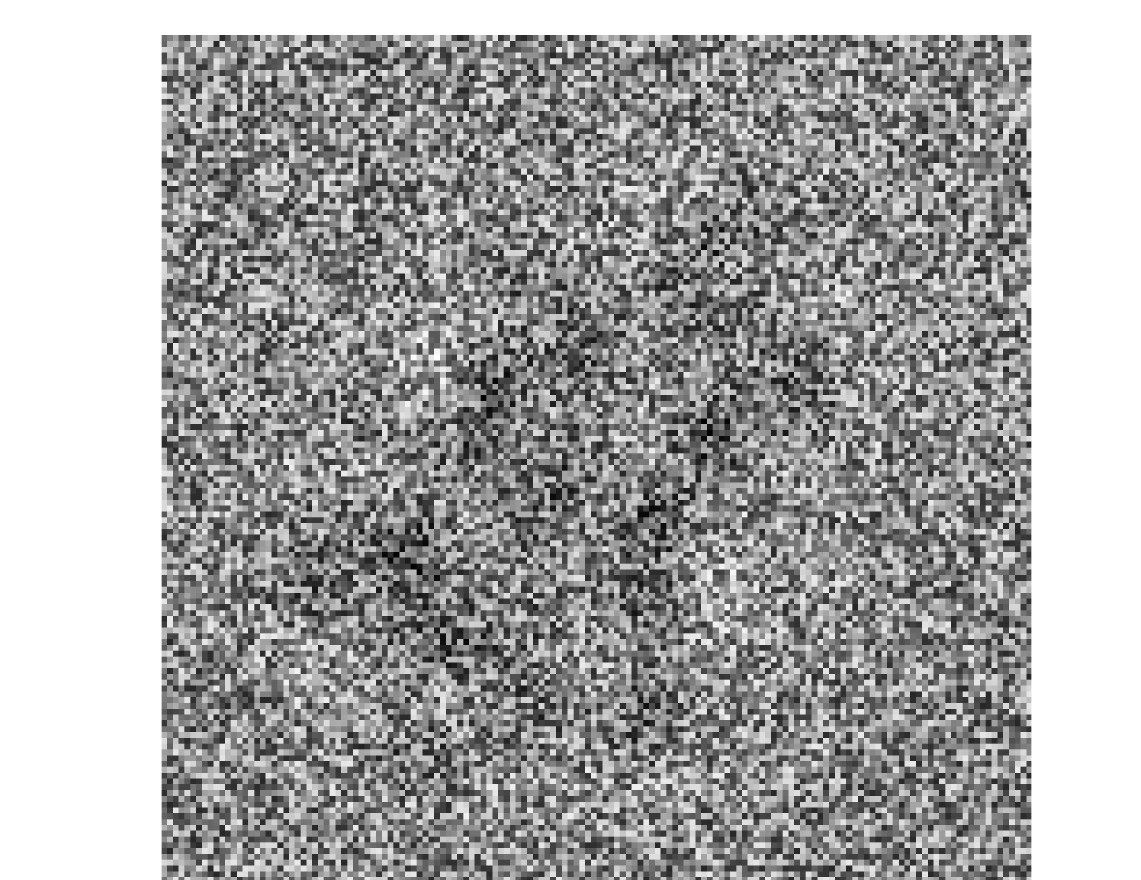}}
\\
(c)\raisebox{-1in}{
\includegraphics[width=.25\linewidth,trim= {6cm 11cm 6cm 11cm},clip]{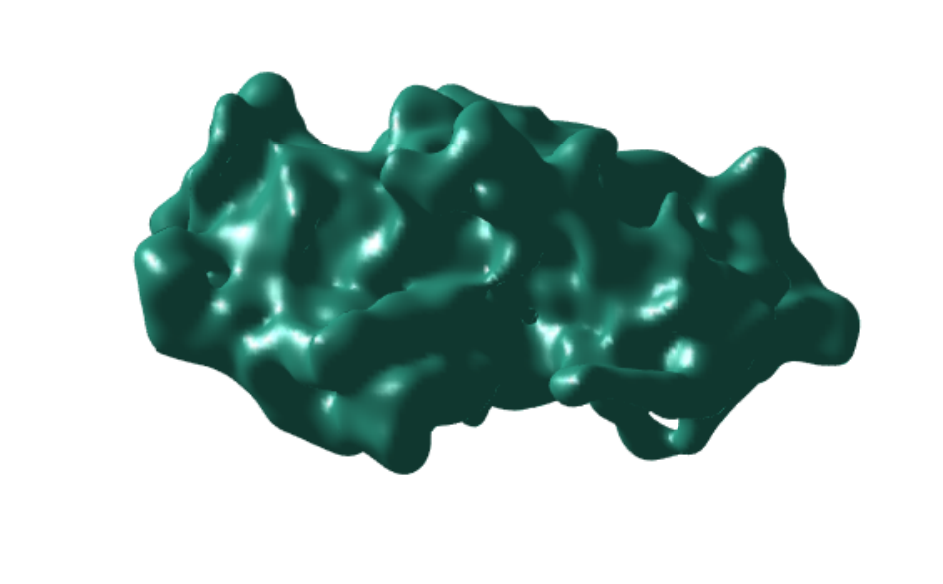}}
&
(d)\raisebox{-1in}{
\includegraphics[width=.25\linewidth,trim= {6cm 11cm 6cm 11cm},clip]{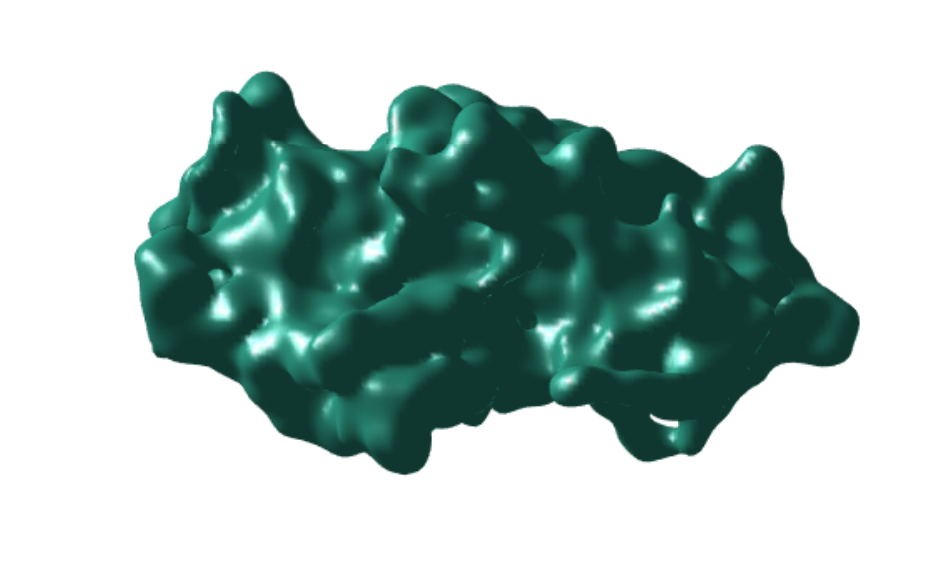}}
\\
(e)\raisebox{-1in}{
\includegraphics[width=.25\linewidth,trim={5cm 9cm 5cm 9cm},clip]{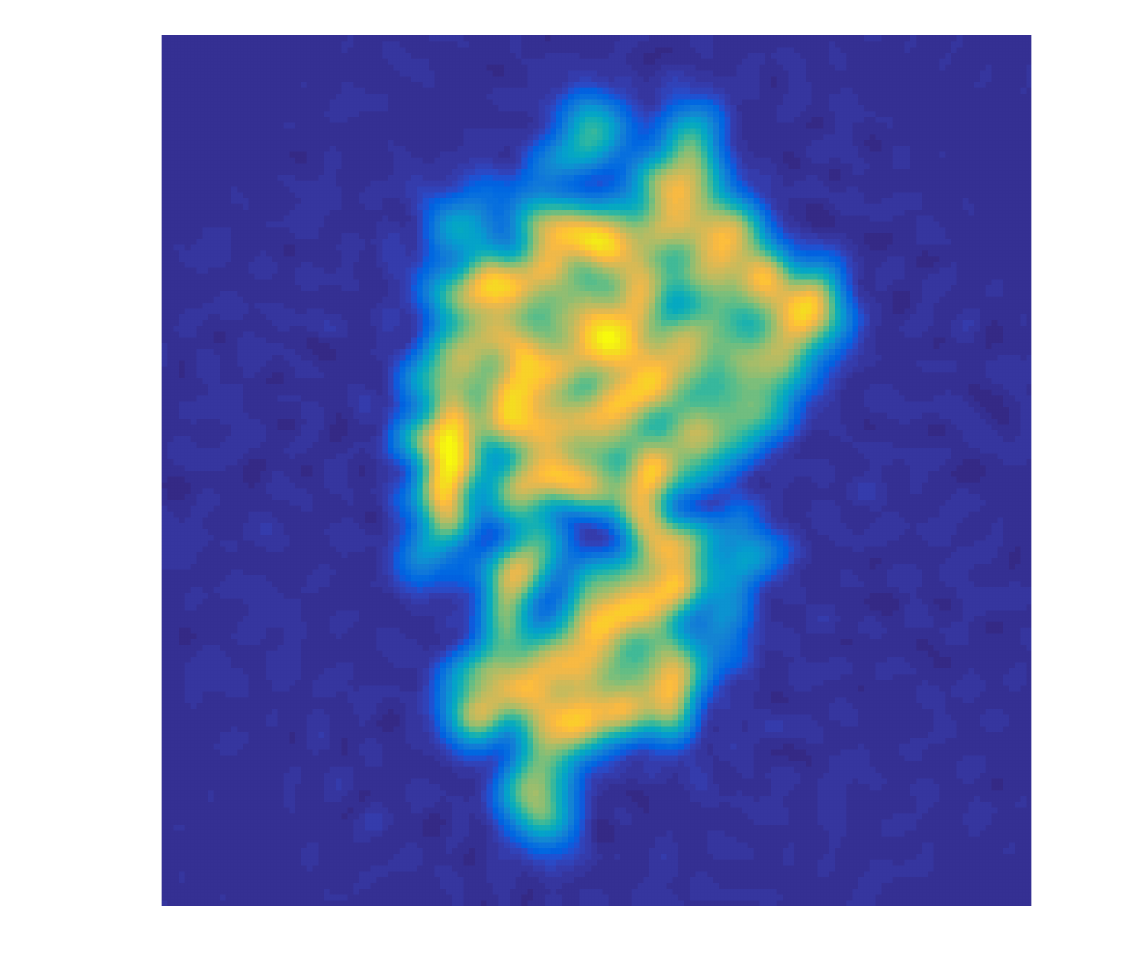}}
&
(f)\raisebox{-1in}{
\includegraphics[width=.25\linewidth,trim={5cm 9cm 5cm 9cm},clip]{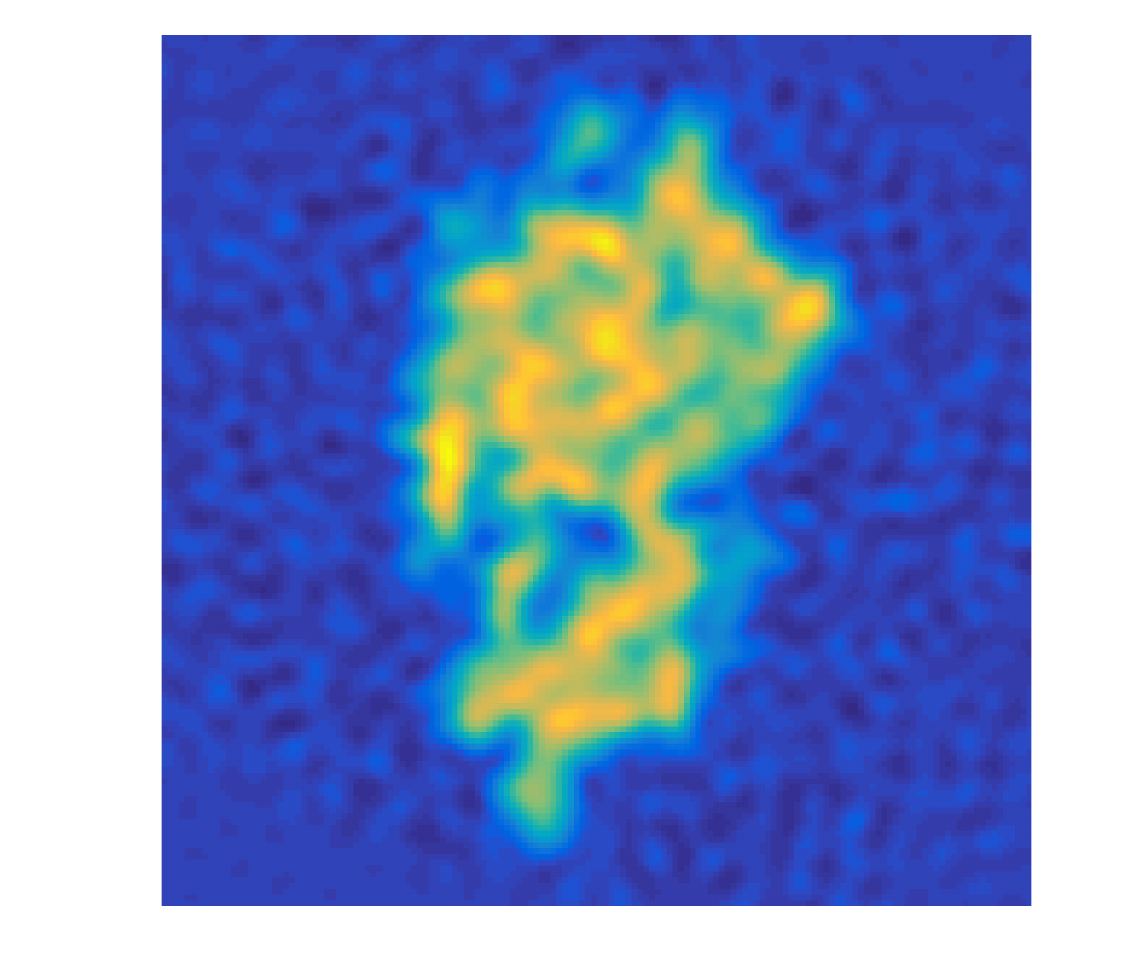}}
\end{tabular}
\ec
\caption{Results for lipoxygenase-1 reconstruction. (a,b): Examples
  of experimental images at SNR 0.5 and 0.05. (c,d):
  Reconstructions by frequency marching using experimental images at SNR 0.5
  and 0.05, shown as isosurfaces. (e,f): Examples of a single slice through the reconstructed $f$.}
\label{figure:lipoxygenease1-results}
\end{figure}

\begin{figure}
\bc
\begin{tabular}{cc}
SNR 0.5 
&
SNR 0.05 
\\
(a)\raisebox{-1in}{
\includegraphics[width=.25\linewidth,trim= {5cm 10cm 5cm 9cm},clip]{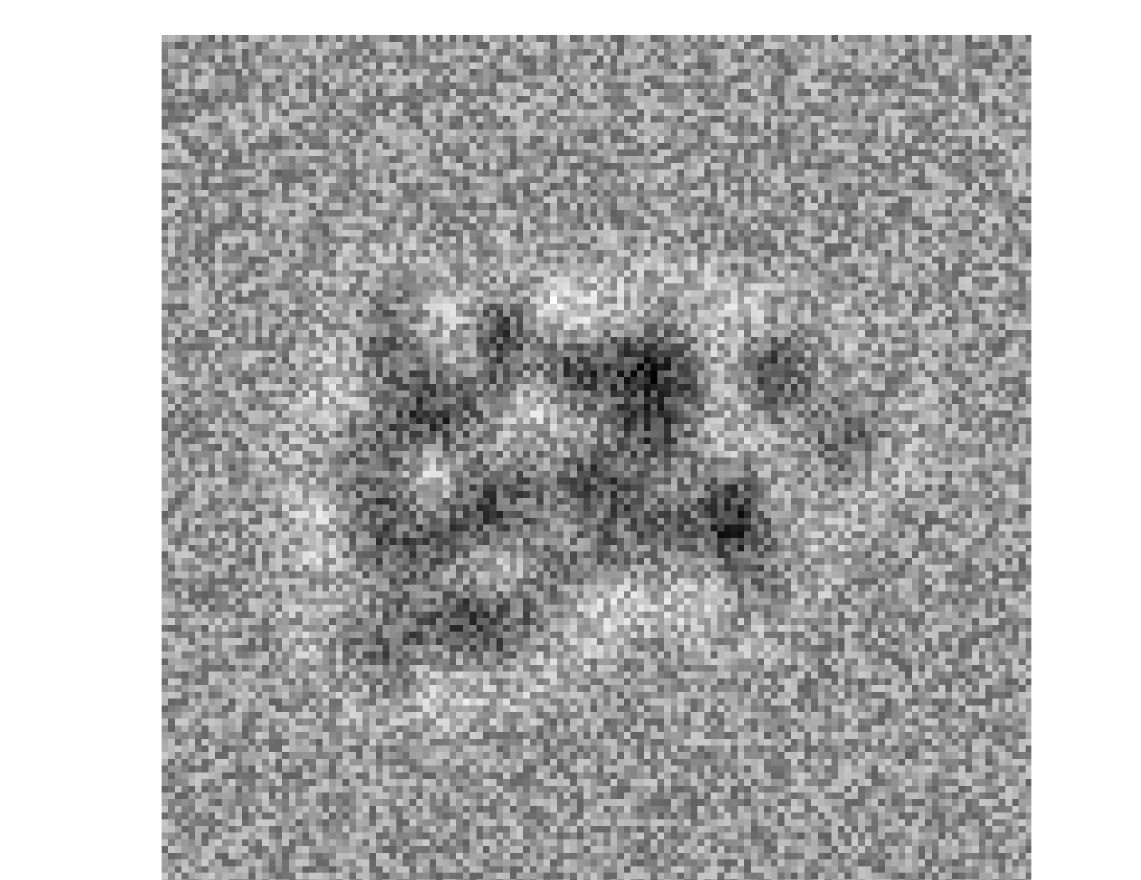}}
&
(b)\raisebox{-1in}{
\includegraphics[width=.25\linewidth,trim= {5cm 10cm 5cm 9cm},clip]{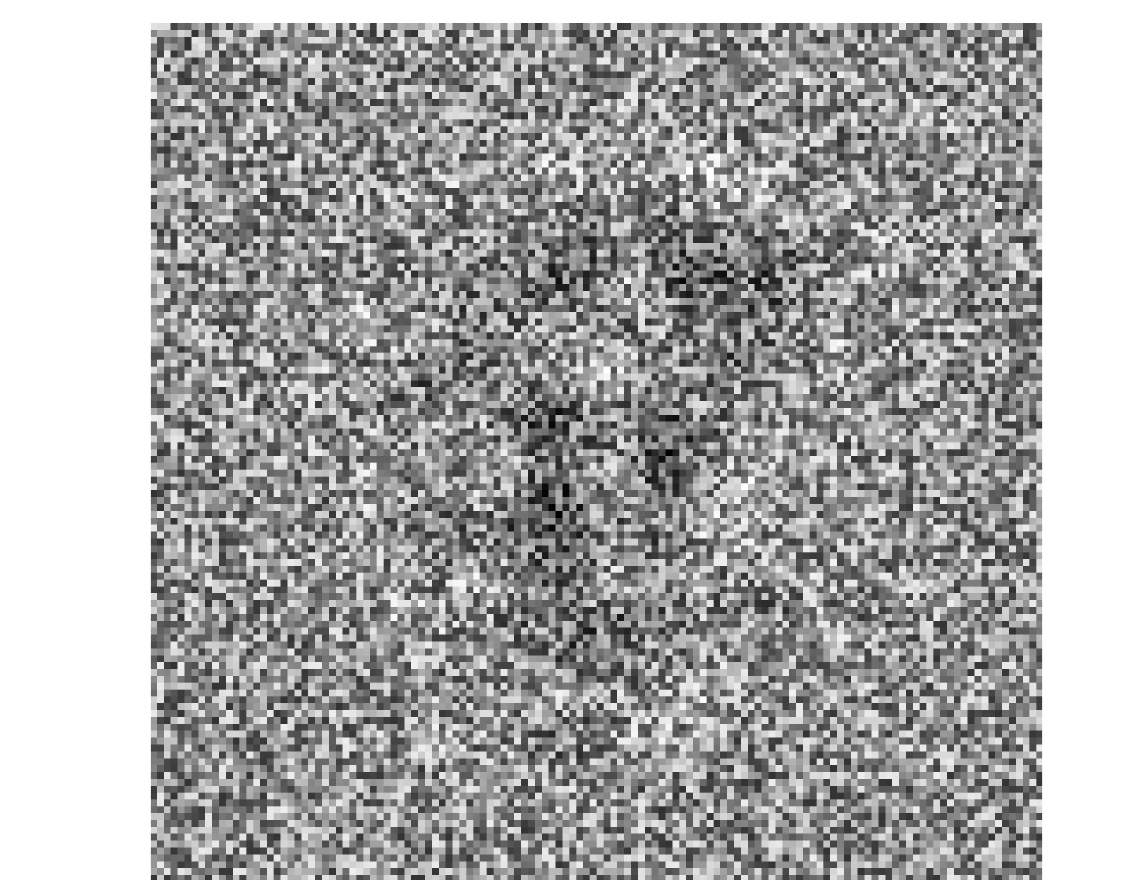}}
\\
(c)\raisebox{-1in}{
\includegraphics[width=.25\linewidth,trim= {6cm 11cm 6cm 10cm},clip]{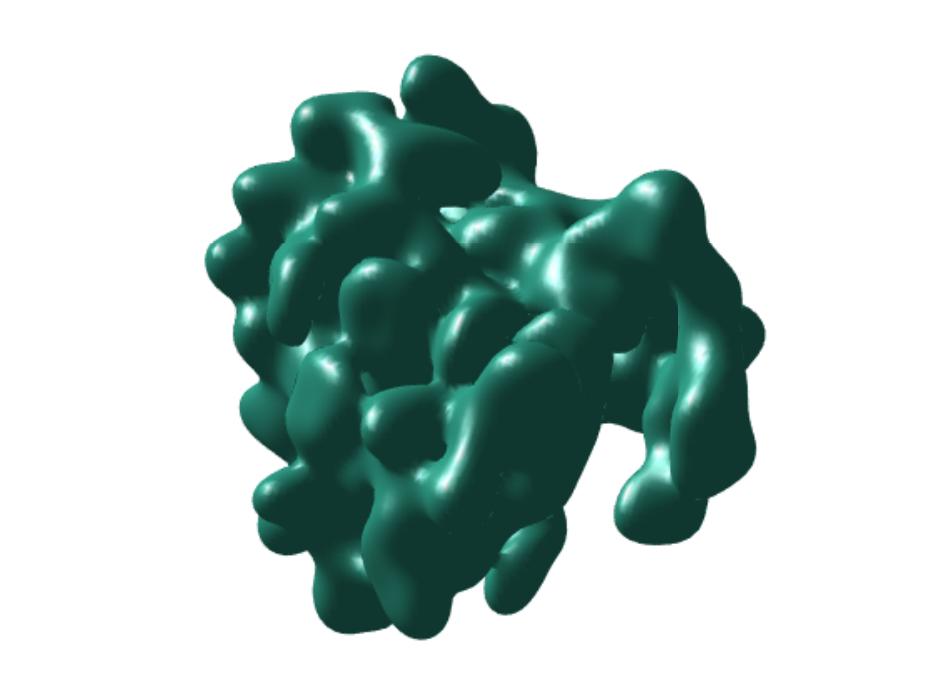}}
&
(d)\raisebox{-1in}{
\includegraphics[width=.25\linewidth,trim= {6cm 11cm 6cm 10cm},clip]{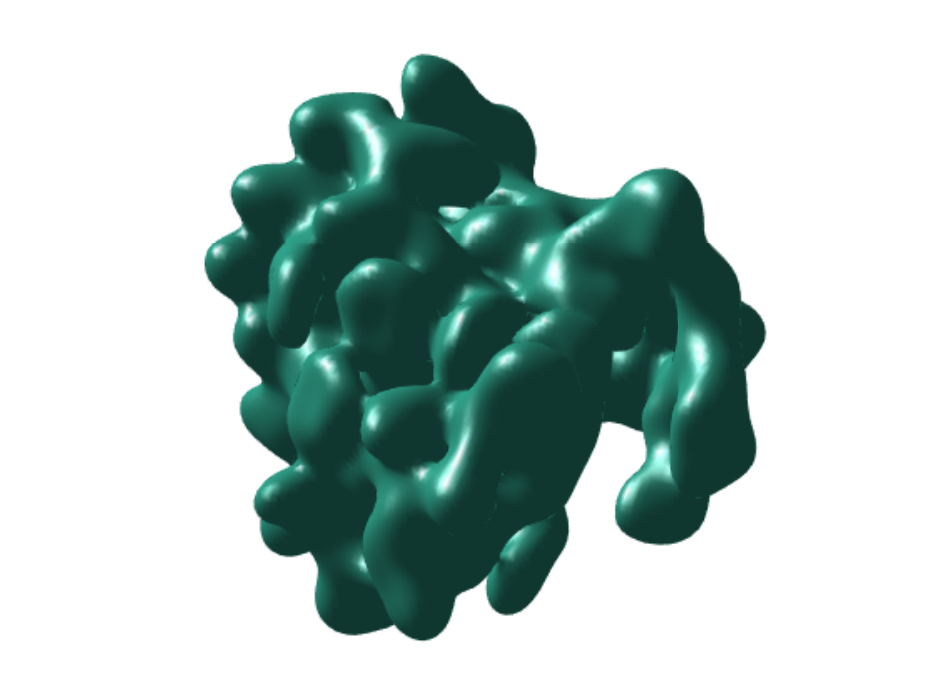}}
\\
(e)\raisebox{-1in}{
\includegraphics[width=.25\linewidth,trim={5cm 9cm 5cm 9cm},clip]{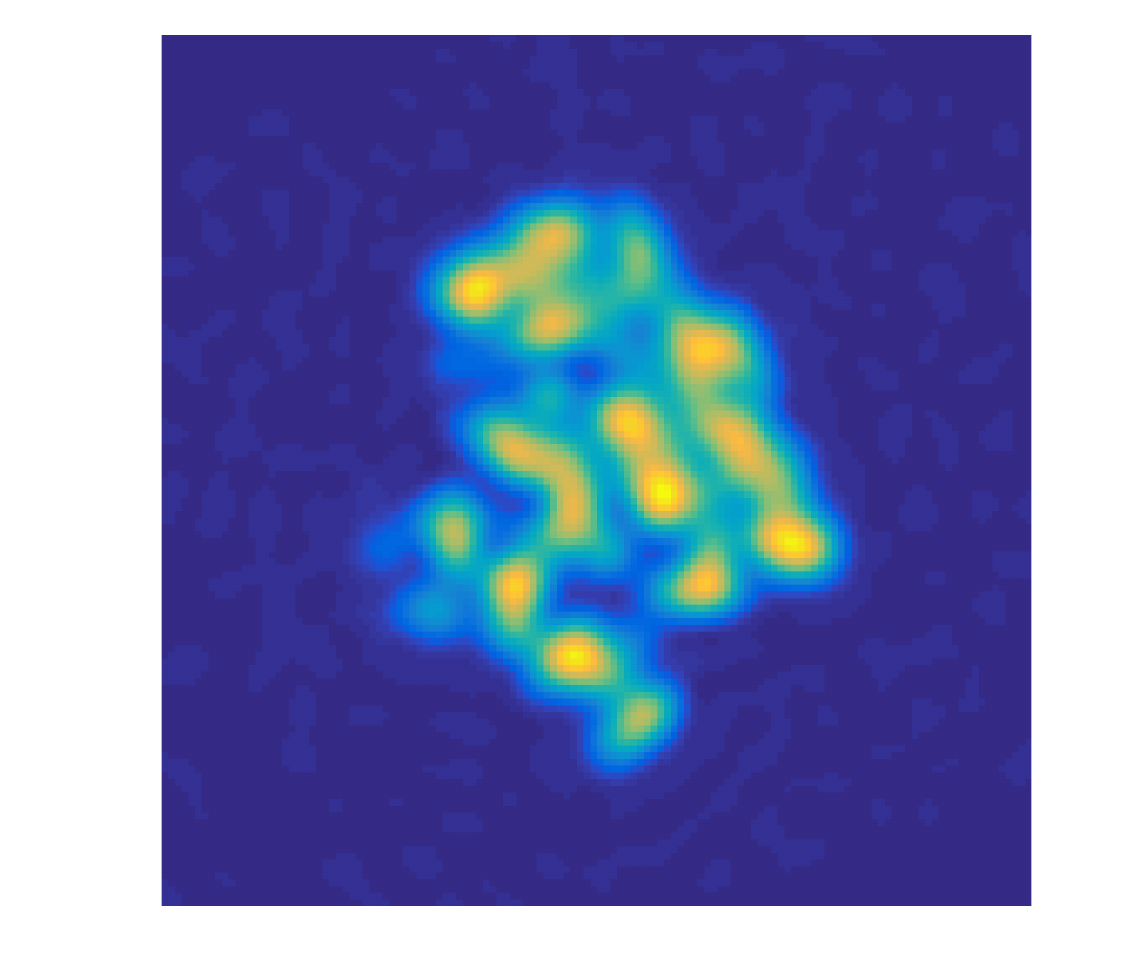}}
&
(f)\raisebox{-1in}{
\includegraphics[width=.25\linewidth,trim={5cm 9cm 5cm 9cm},clip]{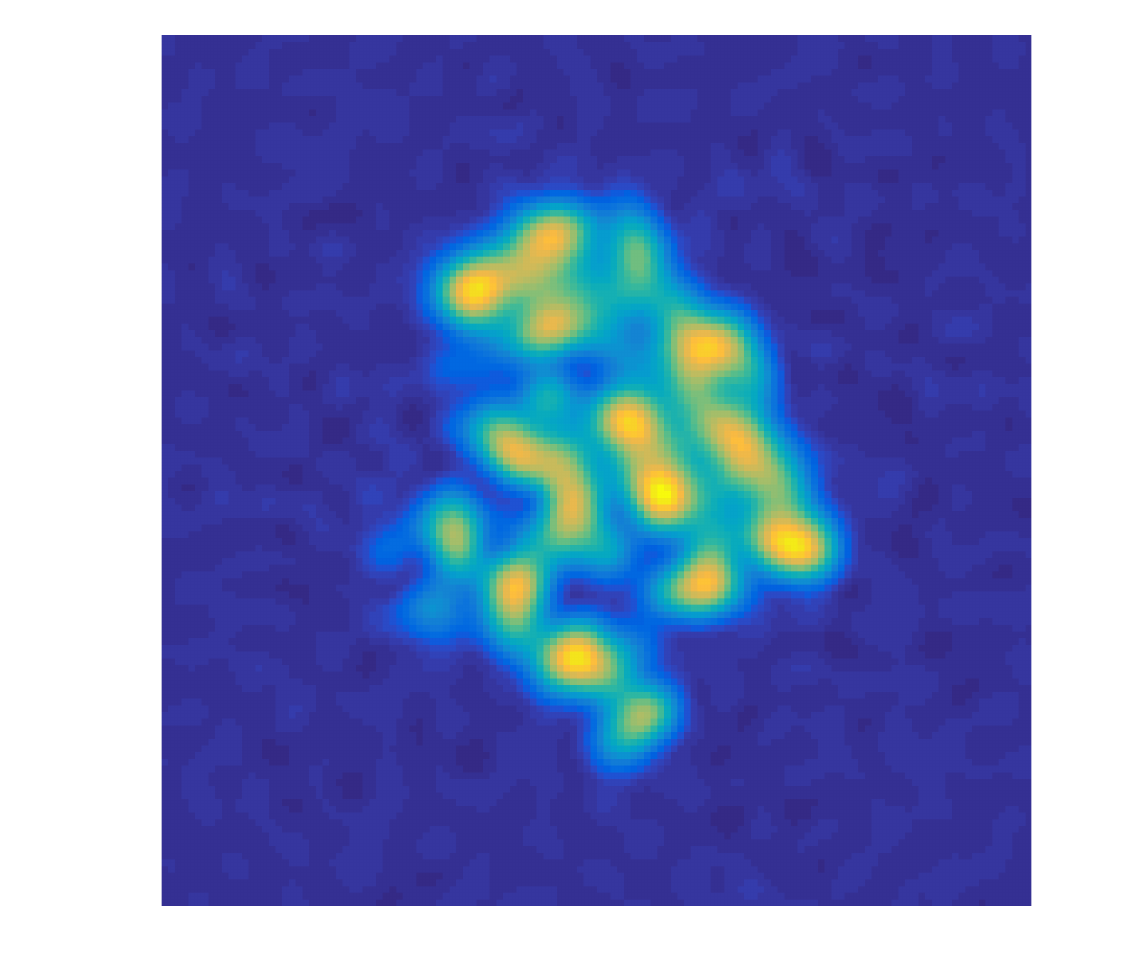}}
\end{tabular}
\ec
\caption{Results for scorpion toxin reconstruction. (a,b): Examples of
  synthetic experimental images at SNR 0.5 and 0.05. 
(c,d): Reconstructions by frequency marching from experimental 
images at SNR 0.5 and 0.05, shown as isosurfaces. 
(e,f): Examples of a single slice through the reconstructed $f$. }
\label{figure:scorpion-toxin-results}
\end{figure}

\begin{figure}  
\includegraphics[width=.33\linewidth]{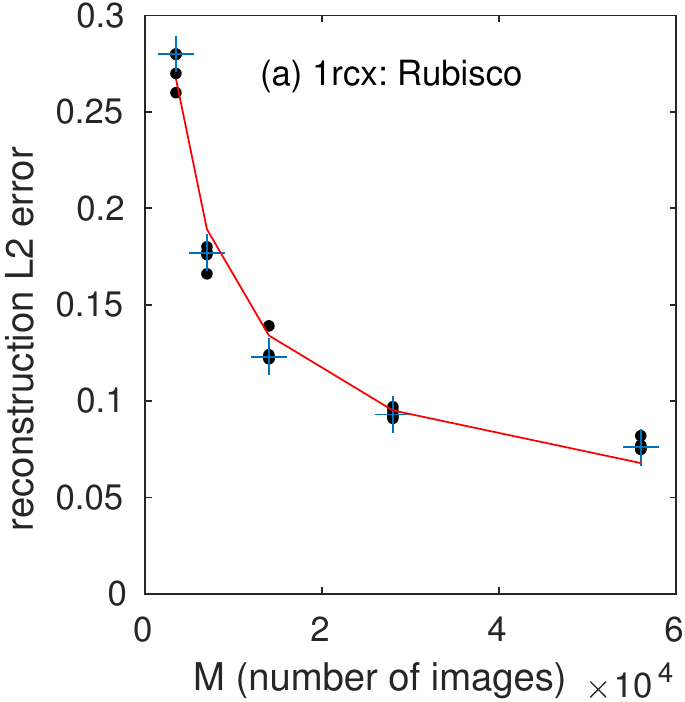}
\includegraphics[width=.33\linewidth]{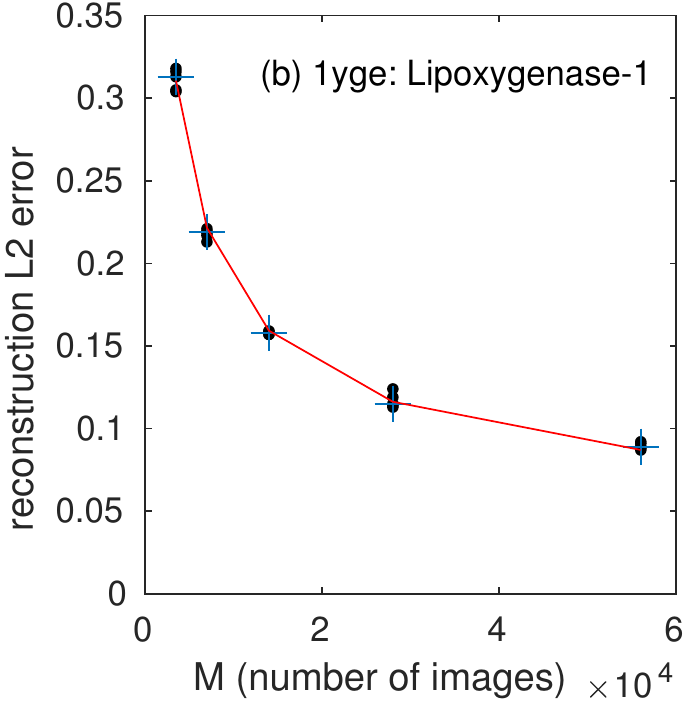}
\includegraphics[width=.33\linewidth]{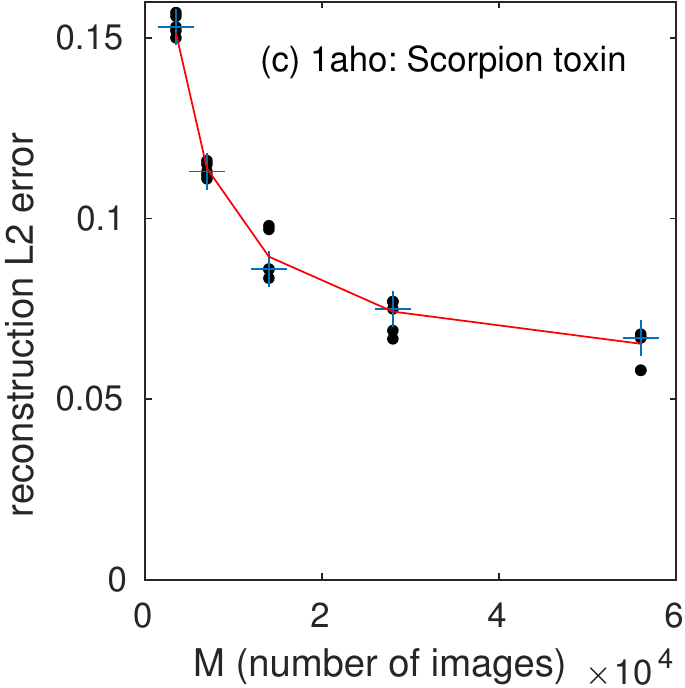}
\caption{Relative $L_2$-error $\epsilon$ of the reconstructed model (see \eqref{err}) vs $M$,
the number of simulated experimental images used for the
reconstruction (black dots), using five runs for each of the three molecules.
The SNR is 0.05.
The error decreases as $M$ increases.
The median is also shown (plus signs),
as is the best fit to the expected form $a_0 + a_1/M$ for the square of the
median errors. (Two of the runs failed with our random initialization
procedure and are not shown on this figure).
\label{figure:l2err-vs-nimages}
}
\end{figure}

Figures \ref{figure:rubisco-results} c) and d),
\ref{figure:lipoxygenease1-results} c) and d) and
\ref{figure:scorpion-toxin-results} c) and d) show the results of the
reconstruction for the three models using experimental images at SNR
0.5 and 0.05 respectively. Out of the 5 runs, we show the
reconstructions with the lowest error. Table \ref{table-timings}
shows that it took approximately 2 hours to do the reconstructions using
frequency marching on 14 cores for all three models. This table also
shows that the time is not significantly affected by the noise levels
in the data. In addition, we expect that the quality of the
reconstruction will improve as we increase the number of experimental
images. {This effect is visible in Figure 
\ref{figure:l2err-vs-nimages},}
where the median squared error of the runs is surprisingly well
fit by the functional form of a constant (accounting for Fourier
and image discretization errors) plus a constant times
$1/M$,
accounting for the usual reduction in statistical variance with
a growing data set size $M$.

\begin{table}
\bc
\begin{tabular}{ll|l|l|l|l|}
&
&
No Noise
&
SNR $0.5$
&
SNR $0.1$
&
SNR $0.05$
\\
\hline
\multirow{2}{*}{Rubisco} 
& 
frequency marching
&
$0.06 \pm 0.001 $
&
$0.06 \pm 0.001 $
&
$0.08 \pm 0.001 $
&
$0.09 \pm 0.004 $
\\
&
known angles
&
$0.05$
&
$0.06$
&
$0.07$
&
$0.08$
\\
\hline
\multirow{2}{*}{Lipoxygenase-1} 
& 
frequency marching
&
$0.05 \pm 0.008$
&
$0.06 \pm 0.002$
&
$0.09 \pm 0.005$
&
$0.12 \pm 0.003$
\\
&
known angles
&
$0.04$
&
$0.05$
&
$0.08$
&
$0.11$
\\
\hline
\multirow{2}{*}{Scorpion toxin} & 
frequency marching
&
$0.05 \pm 0.002$
&
$0.06 \pm 0.003$
&
$0.06 \pm 0.002$
&
$0.07 \pm 0.004$
\\
&
known angles
&
$0.04$
&
$0.04$
&
$0.05$
&
$0.06$
\\
\hline
\end{tabular}
\ec
\vspace{1ex}
\caption{Relative $L_2$ errors $\epsilon$ (see \eqref{err})
of the reconstructions using frequency marching,
compared to the best-possible reconstruction using known angles for the
  experimental images. The frequency marching results are given as
  averages over 5 runs, with estimated standard deviation.}
\label{table-errors}
\end{table}

\begin{table}
\bc
\begin{tabular}{l|l|l|l|l|}
&
No Noise
&
SNR $0.5$
&
SNR $0.1$
&
SNR $0.05$
\\
\hline
Rubisco
&
$1.6 \pm 0.5 $
&
$1.4 \pm 0.3 $
&
$1.2 \pm 0.1 $
&
$1.4 \pm 0.1 $
\\
Lipoxygenase-1
&
$1.1 \pm 0.4 $
&
$1.5 \pm 0.6 $
&
$1.1 \pm 0.4 $
&
$1.2 \pm 0.3 $
\\
Scorpion toxin
&
$0.5 \pm 0.1$
&
$0.5 \pm 0.1$
&
$0.6 \pm 0.1$
&
$0.7 \pm 0.2$
\\
\hline
\end{tabular}
\ec
\vspace{1ex}
\caption{Time in hours for the reconstructions on 14 cores, in the context of
  different levels of noise added to the experimental images. All the
  times are given as averages with estimated standard deviations,
over 5 different runs.}
\label{table-timings}
\end{table}

\section{Conclusions}
\label{conc}

We have presented a fast, robust algorithm for determining 
the three-dimensional structure of proteins using ``single particle" 
cryo-electron microscopy (cryo-EM) data. This problem is typically
formulated in the language of global optimization, leading to a non-convex
objective function for the unknown angles and lateral offsets of each particle
in a collection of noisy projection images.
With hundreds of thousands of images, this leads to a very CPU-intensive 
task. By using a recursive method in the Fourier domain, we have shown
than an essentially deterministic scheme is able to achieve 
high-resolution reconstruction with predictable and modest
computing requirements; in our initial implementation,
fifty thousand images can be processed in about 
one hour on a standard multicore desktop workstation.

One important consequence is that, 
with sufficiently fast reconstruction,
one can imagine making use of multiple runs on the same data,
opening up classical jackknife, bootstrap, or cross-validation statistics 
for quality and resolution assessment. 
In this work, we illustrate the power of doing so with only five
runs.

\vspace{.1in}

\noindent
{\bf Practical considerations}

There are several features which need to be added to the algorithm above
for it to be applicable to experimental data.
Most critically, we need to extend the set of unknowns for each image
to include the translational degrees of freedom (unknown 
lateral offsets for the particle ``center"). 
We believe that frequency marching can again be used to great effect
by defining a fixed search region (either a $3\times 3$ or a $5 \times 5$
offset grid) whose spatial scale corresponds to the current resolution.
That is, we define a spatial step in each lateral direction of the order 
$O(1/k)$, so that large excursions are tested at low resolution and pixel-scale
excursions are tested at high resolution.
This would increase the workload by a constant factor. If done naively,
the increase would correspond to a multiplicative factor of 9 or 25 
on a portion of the code that currently accounts for 50\% of the total time.
We suspect that improvements in the search strategy (using asymptotic 
analysis) can reduce this overhead by a substantial factor and that other 
code optimizations will further reduce the run time (described below).
Thus, we expect that the net increase in cost will be modest over our 
current timings.

\vspace{.1in}

\noindent
{\bf Code optimization}

Significant optimizations are possible to 
further reduce the run time of our template matching and least squares 
reconstruction steps. These include

\begin{enumerate}
\item Using more uniform and sparser sampling in Fourier space:
At frequency $k$, only $O(k^2)$ templates are needed to resolve the 
unknown function $F(k,\theta,\phi)$. Our initial implementation uses
$O(K^2)$ templates where $K$ is the maximum frequency of interest.
A factor of order 5 speedup is available here.
Moreover, by reducing 
the number of azimuthal points near the poles, the number of templates 
can be reduced by an additional factor of $\pi/2$.
\item Updating only the last shell in Fourier space
via the least-squares solve:
Our initial implementation rebuilds the function $F(k,\theta,\phi)$ on
all spheres out to radius $k$ at every iteration 
(see Remark~\ref{r:shell}).
\item Using coarser grids at low resolution in 
frequency marching: At present, we use the same full resolution
spherical grid for for every $k$.
\item Restricting the angle search in ($\alpha,\beta$) to $O(1)$ angles in
the neighborhood of the best guess at the previous frequency $k$:
This would reduce the overall complexity of the step from $O(MK^5)$ to $O(MK^3)$.
Similar strategies have been found to be very effective in
algorithms that marginalize over angles \cite{Brubaker2015a}.
\end{enumerate}

\noindent

While it would increase the cost, it is also worth exploring
the use of marginalization over angles, as in full expectation-maximization
(EM) approaches, as well as the use of sub-grid angle fitting in the projection 
matching step. 
One could also imagine exploring the effectiveness of 
multiple iterations of frequency marching (although preliminary experiments
with noisy data didn't show significant improvements).

Finally, as noted above, our random initialization is subject to 
occasional failure. We will investigate both the possibility of enforcing
convergence by carrying out a more sophisticated nonlinear search for
a low resolution initial guess, or by finding metrics by which to rapidly
discard diverging trajectories.
On a related note, we also plan to develop metrics for discarding images 
that appear to correspond
to erroneously included ``non-particles". In frequency marching, we suspect that
the hallmark of such non-particles (at least in the asymmetric setting)
will be the failure of template matching to find more and more localized
matches. On a more speculative note, we suspect that frequency marching
will be able to handle structural heterogeneity without too many modifications,
at least in the setting where there a finite number of 
well-defined conformations, say $R$,
present in the dataset. The simplest approach would be to assign a conformation
label $(1,\dots,R)$ as an additional parameter for each image 
and to reconstruct $R$ Fourier space models 
in parallel, with random initialization of the labels, and assignments of the
best label  recomputed at each stage.
We will report on all of the above developments at a later date.

\section*{Acknowledgments}
We would like to thank David Hogg, Roy Lederman, Jeremy Magland, 
Christian M\"{u}ller, Adi Rangan, Amit Singer, and Doug Renfrew
for many useful conversations, and the anonymous referees for their
helpful suggestions.

\bibliographystyle{abbrv}
\bibliography{cryo_bib_copy}

\begin{thebibliography}{10}

\bibitem{BaoLi2015}
G.~Bao, P.~Li, J.~Lin, and F.~Triki.
\newblock Inverse scattering problems with multi-frequencies.
\newblock {\em Inverse Problems}, 31(9):093001, 2015.

\bibitem{Bell2016}
J.~M. Bell, M.~Chen, P.~R. Baldwin, and S.~J. Ludtke.
\newblock {High resolution single particle refinement in EMAN2.1}.
\newblock {\em Methods}, 100:25--34, 2016.

\bibitem{pdb}
H.~Berman, J.~Westbrook, Z.~Feng, G.~Gilliland, T.~N. Bhat, H.~Weissig, I.~N.
  Shindyalov, and P.~E. Bourne.
\newblock The protein data bank.
\newblock {\em Nucleic Acids Research}, 28:235--242, 2000.
\newblock {\tt http://www.rcsb.org}.

\bibitem{Borges2016}
C.~Borges, A.~Gillman, and L.~Greengard.
\newblock High resolution inverse scattering in two dimensions using recursive
  linearization.
\newblock {\em SIAM J. Imaging Sci.}, 10(2):641--664, 2017.

\bibitem{Brubaker2015a}
M.~A. Brubaker, A.~Punjani, and D.~J. Fleet.
\newblock Building proteins in a day: Efficient {3D} molecular reconstruction.
\newblock In {\em 2015 IEEE Conference on Computer Vision and Pattern
  Recognition (CVPR)}, pages 3099--3108, June 2015.

\bibitem{Chenrep1088}
Y.~Chen.
\newblock Recursive linearization for inverse scattering.
\newblock Technical Report Yale Research Report/DCS/RR-1088, Department of
  Computer Science, Yale University, New Haven, CT, October 1995.

\bibitem{Chen}
Y.~Chen.
\newblock Inverse scattering via {H}eisenberg's uncertainty principle.
\newblock {\em Inverse Problems}, 13:253--282, 1997.

\bibitem{Cheng2015}
Y.~Cheng, N.~Grigorieff, P.~A. Penczek, and T.~Walz.
\newblock A primer to single-particle cryo-electron microscopy.
\newblock {\em Cell}, 161:439--449, 2015.

\bibitem{Coifman2010}
R.~Coifman, Y.~Shkolnisky, F.~J. Sigworth, and A.~Singer.
\newblock Reference free structure determination through eigenvectors of center
  of mass operators.
\newblock {\em Appl. Comput. Harmonic Anal.}, 47(3):296--312, 2010.

\bibitem{Cong2003}
Y.~Cong, J.~A. Kovacs, and W.~Wriggers.
\newblock {2D fast rotational matching for image processing of biophysical
  data}.
\newblock {\em Journal of Structural Biology}, 144(1-2):51--60, 2003.

\bibitem{Dvornek2015}
N.~C. Dvornek, F.~J. Sigworth, and H.~D. Tagare.
\newblock {SubspaceEM: A fast maximum-a-posteriori algorithm for cryo-EM single
  particle reconstruction}.
\newblock {\em Journal of Structural Biology}, 190:200--214, 2015.

\bibitem{Elmlund2012}
D.~Elmlund and H.~Elmlund.
\newblock {SIMPLE: Software for ab initio reconstruction of heterogeneous
  single-particles}.
\newblock {\em Journal of Structural Biology}, 180:420--427, 2012.

\bibitem{SPIDER}
J.~Frank, B.~Shimkin, and H.~Dowse.
\newblock {SPIDER}---a modular software system for electron image processing.
\newblock {\em Ultramicroscopy}, 6:343--358, 1981.

\bibitem{Goncharov1987}
A.~B. Goncharov.
\newblock {Methods of integral geometry and finding the relative orientation of
  identical particles arbitrarily arranged in a plane from their projections
  onto a stright line}.
\newblock {\em Dokl. Akad. Nauk SSSR}, 293:355--58, 1987.

\bibitem{Goncharov1988}
A.~B. Goncharov and M.~S. Gelfand.
\newblock {Determination of mutual orientation of identical particles from
  their projections by the moments method}.
\newblock {\em Ultramicroscopy}, 25:317--28, 1988.

\bibitem{nufft}
L.~Greengard and J.-Y. Lee.
\newblock Accelerating the nonuniform fast {F}ourier transform.
\newblock {\em SIAM Review}, 46(3):443--454, 2004.

\bibitem{Grigorieff2007}
N.~Grigorieff.
\newblock {FREALIGN: High-resolution refinement of single particle structures}.
\newblock {\em Journal of Structural Biology}, 157(1):117--125, 2007.

\bibitem{Grigorieff2016}
N.~Grigorieff.
\newblock Frealign: An exploratory tool for single-particle {cryo-EM}.
\newblock In R.~A. Crowther, editor, {\em The Resolution Revolution: Recent
  Advances In {cryoEM}}, volume 579 of {\em Methods in Enzymology}, pages
  191--226. Academic Press, 2016.

\bibitem{Henderson2013}
R.~Henderson.
\newblock {Avoiding the pitfalls of single particle cryo-electron microscopy:
  Einstein from noise.}
\newblock {\em Proceedings of the National Academy of Sciences of the United
  States of America}, 110(45):18037--41, 2013.

\bibitem{Henderson2012}
R.~Henderson, A.~Sali, M.~L. Baker, B.~Carragher, B.~Devkota, K.~H. Downing,
  E.~H. Egelman, Z.~Feng, J.~Frank, N.~Grigorieff, W.~Jiang, S.~J. Ludtke,
  O.~Medalia, P.~A. Penczek, P.~B. Rosenthal, M.~G. Rossmann, M.~F. Schmid,
  G.~F. Schr{\"{o}}der, A.~C. Steven, D.~L. Stokes, J.~D. Westbrook,
  W.~Wriggers, H.~Yang, J.~Young, H.~M. Berman, W.~Chiu, G.~J. Kleywegt, and
  C.~L. Lawson.
\newblock {Outcome of the first electron microscopy validation task force
  meeting}.
\newblock {\em Structure}, 20(2):205--214, 2012.

\bibitem{Heymann2015}
J.~B. Heymann.
\newblock {Validation of 3D EM Reconstructions: The Phantom in the Noise.}
\newblock {\em AIMS biophysics}, 2:21--35, 2015.

\bibitem{SPARX}
M.~Hohn, G.~Tang, G.~Goodyear, P.~Baldwin, Z.~Huang, P.~Penczek, C.~Yang,
  R.~Glaeser, P.~Adams, and S.~Ludtke.
\newblock Sparx, a new environment for cryo-em image processing.
\newblock {\em J. Struct. Biol.}, 157:47--55, 2007.

\bibitem{Joubert2015}
P.~Joubert and M.~Habeck.
\newblock {B}ayesian inference of initial models in cryo-electron microscopy
  using pseudo-atoms.
\newblock {\em Biophysical Journal}, 108:1165--1175, 2015.

\bibitem{Joyeux2002}
L.~Joyeux and P.~A. Penczek.
\newblock {Efficiency of 2D alignment methods}.
\newblock {\em Ultramicroscopy}, 92(2):33--46, 2002.

\bibitem{Kimanius2016}
D.~Kimanius, B.~O. Forsberg, S.~H. Scheres, and E.~Lindehl.
\newblock {Accelerated cryo-EM structure determination with parallelisation
  using GPUs in RELION-2}.
\newblock {\em eLife}, 5:e18722, 2016.

\bibitem{Lyumkis2013}
D.~Lyumkis, A.~F. Brilot, D.~L. Theobald, and N.~Grigorieff.
\newblock {Likelihood-based classification of cryo-EM images using FREALIGN}.
\newblock {\em Journal of Structural Biology}, 183(3):377--388, 2013.

\bibitem{Milne2013}
J.~L.~S. Milne, M.~J. Borgnia, A.~Bartesaghi, E.~E.~H. Tran, L.~A. Earl, D.~M.
  Schauder, J.~Lengyel, J.~Pierson, A.~Patwardhan, and S.~Subramaniam.
\newblock {Cryo-electron microscopy - A primer for the non-microscopist}.
\newblock {\em FEBS Journal}, 280:28--45, 2013.

\bibitem{Mindell2003}
J.~A. Mindell and N.~Grigorieff.
\newblock {Accurate determination of local defocus and specimen tilt in
  electron microscopy}.
\newblock {\em Journal of Structural Biology}, 142:334--347, 2003.

\bibitem{Natterer}
F.~Natterer.
\newblock {\em The Mathematics of Computerized Tomography}.
\newblock SIAM, 2001.

\bibitem{Nogales2016}
E.~Nogales.
\newblock {The development of cryo-EM into a mainstream structural biology
  technique}.
\newblock {\em Nat. Meth.}, 13:24--27, 2016.

\bibitem{Penczek1992}
P.~Penczek, M.~Radermacher, and J.~Frank.
\newblock {Three-dimensional reconstruction of single particles embedded in
  ice}.
\newblock {\em Ultramicroscopy}, 40:33--53, 1992.

\bibitem{Frank2006}
P.~A. Penczek.
\newblock {\em Three-Dimensional Electron Microscopy of Macromolecular
  Assemblies: Visualization of Biological Molecules in Their Native State}.
\newblock Oxford University Press, 2006.

\bibitem{Penczek2010rev}
P.~A. Penczek.
\newblock Resolution measures in molecular electron microscopy.
\newblock {\em Methods Enzymol.}, 482:73--100, 2010.

\bibitem{Penczek1994}
P.~A. Penczek, R.~A. Grassucci, and J.~Frank.
\newblock {The ribosome at improved resolution: New techniques for merging and
  orientation refinement in 3D cryo-electron microscopy of biological
  particles}.
\newblock {\em Ultramicroscopy}, 53(3):251--270, 1994.

\bibitem{Penczek1996}
P.~A. Penczek, J.~Zhu, and J.~Frank.
\newblock A common-lines based method for determining orientations for {$N >
  3$} particle projections simultaneously.
\newblock {\em Ultramicroscopy}, 63:205--218, 1996.

\bibitem{Provencher1988}
S.~W. Provencher and R.~H. Vogel.
\newblock {Three-dimensional reconstruction from electron micrographs of
  disordered specimens I. Method}.
\newblock {\em Ultramicroscopy}, 25:209--221, 1988.

\bibitem{Rosenthal2015}
P.~B. Rosenthal and J.~L. Rubinstein.
\newblock {Validating maps from single particle electron cryomicroscopy}.
\newblock {\em Current Opinion in Structural Biology}, 34:135--144, 2015.

\bibitem{Belnap2010}
E.~Sanz-Garc{\'{i}}a, A.~B. Stewart, and D.~M. Belnap.
\newblock {The random-model method enables ab initio three-dimensional
  reconstruction of asymmetric particles and determination of particle
  symmetry}.
\newblock {\em Journal of Structural Biology}, 171(2):216--222, 2010.

\bibitem{Scheres2012}
S.~H.~W. Scheres.
\newblock A {B}ayesian view on {cryo-EM} structure determination.
\newblock {\em Journal of Molecular Biology}, 415:406--418, 2012.

\bibitem{Scheres2012b}
S.~H.~W. Scheres.
\newblock {RELION: Implementation of a Bayesian approach to cryo-EM structure
  determination}.
\newblock {\em Journal of Structural Biology}, 180(3):519--530, 2012.

\bibitem{Scheres2009}
S.~H.~W. Scheres, M.~Valle, P.~Grob, E.~Nogales, and J.-M. Carazo.
\newblock Maximum likelihood refinement of electron microscopy data with
  normalization errors.
\newblock {\em J.\ Struct.\ Biol.}, 166(2):234--240, 2009.

\bibitem{Shkolnisky2012}
Y.~Shkolnisky and A.~Singer.
\newblock {Viewing direction estimation in cryo-EM using synchronization.}
\newblock {\em SIAM J. Imaging Sci.}, 5(3):1088--1110, 2012.

\bibitem{Sigworth1998}
F.~J. Sigworth.
\newblock {A maximum-likelihood approach to single-particle image refinement}.
\newblock {\em Journal of structural biology}, 122(3):328--39, 1998.

\bibitem{Sigworth2010}
F.~J. Sigworth, D.~P.C., J.-M. Carazo, and S.~Scheres.
\newblock {An introduction to maximum-likelihood methods in Cryo-EM}.
\newblock In {\em Methods in Enzymology. Cryo-EM, Part B: 3D reconstruction},
  pages 263--94. Academic Press., 2010.

\bibitem{Singer2009}
A.~Singer, R.~R. Coifman, F.~J. Sigworth, D.~W. Chester, and Y.~Shkolnisky.
\newblock Detecting consistent common lines in {cryo-EM} by voting.
\newblock {\em J. Struct. Biol.}, 169:312--322, 2009.

\bibitem{Singer2011}
A.~Singer and Y.~Shkolnisky.
\newblock Three-dimensional structure determination from common lines in
  {cryo-EM} by eigenvectors and semidefinite programming.
\newblock {\em SIAM J. Imaging Sci.}, 4(2):543--72, 2011.

\bibitem{EMAN2}
G.~Tang, L.~Peng, P.~Baldwin, D.~Mann, W.~Jiang, I.~Rees, and S.~Ludtke.
\newblock {EMAN2}: an extensible image processing suite for electron
  microscopy.
\newblock {\em J. Struct. Biol.}, 157:38--46, 2007.

\bibitem{PTRtref}
L.~N. Trefethen and J.~A.~C. Weideman.
\newblock The exponentially convergent trapezoidal rule.
\newblock {\em SIAM Review}, 56(3):385--458, 2014.

\bibitem{Vainshtein1986}
B.~Vainshtein and A.~Goncharov.
\newblock {Determination of the spatial orientation of arbitrary arranged
  identical particles of an unknown structure from their projections}.
\newblock In {\em Proceedings of the 11th International Congress on Electron
  Microscopy}, pages 459--60, 1986.

\bibitem{vanHeel1987}
M.~van Heel.
\newblock {Angular reconstruction: a posteriori assignment of projection
  directions for 3D reconstruction}.
\newblock {\em Ultramicroscopy}, 21:111--23, 1987.

\bibitem{VanHeel1997}
M.~van Heel, E.~V. Orlova, G.~Harauz, H.~Stark, P.~Dube, F.~Zemlin, and
  M.~Schatz.
\newblock {Angular Reconstitution in Three-Dimensional Electron Microscopy:
  Historical and Theoretical Aspects}.
\newblock {\em Scanning Microscopy}, 11:195--210, 1997.

\bibitem{Vinothkumar2016}
K.~R. Vinothkumar and R.~Henderson.
\newblock {Single particle electron cryomicroscopy: trends, issues and future
  perspective}.
\newblock {\em Quarterly Reviews of Biophysics}, 49:e13, 2016.

\bibitem{Vogel1988}
R.~H. Vogel and S.~W. Provencher.
\newblock {Three-dimensional reconstruction from electron micrographs of
  disordered specimens II. Implementation and results}.
\newblock {\em Ultramicroscopy}, 25:223--240, 1988.

\bibitem{wade}
R.~H. Wade.
\newblock A brief look at imaging and contrast transfer.
\newblock {\em Ultramicroscopy}, 46:145--156, 1992.

\bibitem{Wang2013}
L.~Wang, A.~Singer, and Z.~Wen.
\newblock {Orientation determination from cryo-EM images using least unsquared
  deviations}.
\newblock {\em SIAM J. Imaging Sci.}, 6(4):2450--83, 2013.

\bibitem{kink}
J.~A.~C. Weideman and L.~N. Trefethen.
\newblock The kink phenomenon in {Fej\'er} and {Clenshaw--Curtis} quadrature.
\newblock {\em Numer.\ Math.}, 107:707--727, 2007.

\end{thebibliography}

\end{document}